\RequirePackage[l2tabu, orthodox]{nag}        
\documentclass[10pt]{article}

\usepackage{graphicx,caption}
\usepackage{booktabs}
\usepackage[T1]{fontenc}
\usepackage{lmodern}
\usepackage[numbers]{natbib}
\usepackage{microtype}
\usepackage{systeme}
\usepackage{ulem}

\usepackage{sistyle}
\SIthousandsep{,}
\usepackage[frenchmath]{mathastext}                    
\usepackage{amsmath}                        
\numberwithin{equation}{section}
\usepackage{amsfonts}
\usepackage{bbold}
\usepackage{commath}
\usepackage{bm}
\usepackage{stackrel}
\usepackage{graphicx}                      
\usepackage{graphicx}
\usepackage{caption}
\usepackage[skip=0.5ex]{subcaption}
\usepackage[plainpages=false, pdfpagelabels]{hyperref} 
    \hypersetup{
        colorlinks   = true,
        citecolor    = RoyalBlue, 
        linkcolor    = RubineRed, 
        urlcolor     = Turquoise
    }
\usepackage[dvipsnames]{xcolor}
\usepackage[all]{xy}
\usepackage{algorithm}                      
\usepackage{algpseudocode}
\usepackage{amssymb}                        
\usepackage[mathscr]{eucal}                 
\usepackage{dsfont}                                                    
\usepackage[paperwidth=8.5in,paperheight=11.00in,top=1.25in, bottom=1.25in, left=1.00in, right=1.00in]{geometry}
\usepackage{mathtools}                      
\mathtoolsset{showonlyrefs=true}            
\usepackage{fixltx2e,amsmath}               
\MakeRobust{\eqref}
\usepackage{microtype}                                            
\usepackage{amsthm}                         
\allowdisplaybreaks                         
\theoremstyle{plain}
\newtheorem{theorem}{Theorem}
\numberwithin{theorem}{section}

\theoremstyle{definition}

\newtheorem{remark}[theorem]{Remark}

\newcommand{\X}{\mathcal{X}}

\newcommand{\A}{\mathcal{A}}

\newcommand{\E}{\mathds{E}}

\author{
Andrea Angiuli\thanks{Prime Machine Learning Team, Amazon. 320 Westlake Ave N, SEA83, Seattle, WA, 98109 (E-mail: \href{mailto:aangiuli@amazon.com}{aangiuli@amazon.com}). The work presented here does not relate to this author's position at Amazon. }
\and Nils Detering\thanks{Department of Statistics and Applied Probability, South Hall, University of California, Santa Barbara, CA 93106, USA (E-mail: \href{mailto:detering@pstat.ucsb.edu}{detering@pstat.ucsb.edu}). } 
  \and Jean-Pierre Fouque\thanks{Department of Statistics and Applied Probability, South Hall, University of California, Santa Barbara, CA 93106, USA (E-mail: \href{mailto:fouque@pstat.ucsb.edu}{fouque@pstat.ucsb.edu}). Work supported by NSF grants DMS-1814091 and DMS-1953035.} 
  \and Mathieu Lauri\`ere\thanks{Mathematics and Data Science, NYU Shanghai, China (E-mail: 
  \href{mailto:mathieu.lauriere@gmail.com}{mathieu.lauriere@nyu.edu}).}
  \and  Jimin Lin\thanks{Department of Statistics and Applied Probability, South Hall, University of California, Santa Barbara, CA 93106, USA (E-mail: \href{mailto:jiminlin@pstat.ucsb.edu}{jiminlin@pstat.ucsb.edu}). Use was made of computational facilities purchased with funds from the National Science Foundation (CNS-1725797) and administered by the Center for Scientific Computing (CSC). The CSC is supported by the California NanoSystems Institute and the Materials Research Science and Engineering Center (MRSEC; NSF DMR 1720256) at UC Santa Barbara.} 
  }
\begin{document}
\title{Reinforcement Learning Algorithm for Mixed Mean Field Control Games}

\maketitle

\begin{abstract}
We present a new combined \textit{mean field control game} (MFCG) problem which can be interpreted as a competitive game between collaborating groups and its solution as a Nash equilibrium between groups. Players coordinate their strategies within each group. An example is a modification of the classical trader's problem. Groups of traders maximize their wealth. They face cost for their transactions, for their own terminal positions, and for the average holding within their group. The asset price is impacted by the trades of all agents. We propose a three-timescale reinforcement learning algorithm to approximate the solution of such MFCG problems. We test the algorithm on benchmark linear-quadratic specifications for which we provide analytic solutions. 
\end{abstract}

\section{Introduction}
Mean field approaches are based on the idea that the main properties of large coupled systems of entities (e.g. agents, players, or particles) can be described by the distribution of one representative entity. To answer many questions related to the system, it is not required to know the individual states of all entities but only the distribution of their representative. This reduces significantly the complexity of large systems. 

Mean field approaches were first introduced in the context of statistical physics where propagation of chaos among particles was studied. Under mild assumptions, in a system of particles described by a large system of diffusion processes, the location of one particle becomes independent of the others as the size of the system grows \cite{10.1007/BFb0085169}. In the following we think of the entities of the system being agents or players and we have mainly financial applications in mind.  
Mean field ideas have later been adapted to differential games with large number of agents in the cooperative setting (mean field control, MFC), and in the competitive framework (mean field games, MFG) \cite{MR2295621,carmona2018probabilisticI-II,MR3134900}. MFC and MFG problems arise in a number of applications ranging from engineering to economics. Mean field type games (MFTG)~\cite{tembine} are games with a finite number of players who are of mean field type, \textit{i.e.}, their dynamics and cost functions may depend on their own distribution. 

Recently numerical solution of MFC and MFG problems has received greater attention \cite{MR2679575,angiuli2019cemracs,MR3258261_DPMFC_CRAS,MR3608094,Han-Hu-2020,guo2019learning,yang2018deep}; see \textit{e.g.}~\cite{achdoulauriere2020mfgnumerical} for a survey. Classical methods of optimization theory have been complemented by deep neural networks \cite{MR2346927,MR2352434,carmona2019convergence-I,carmona2019convergence-II,fouque2019deep} and by Reinforcement Learning (RL) approaches which aim at calculating optimal strategies without the precise knowledge of the underlying model \cite{SubramanianMahajan-2018-RLstatioMFG,yang2018mean,CarmonaLauriereTan-2019-MQFRL,CarmonaLauriereTan-2019-LQMFRL,mguni2018decentralised,elie2020convergence}. 

In \cite{AFL2021}, a unified reinforcement Q-learning algorithm is proposed to solve MFG and MFC problems based on the ratio of two learning rates, one for the decision Q-matrix and the other for the distribution of the population. In the present paper, we argue that this algorithm can be adapted for solving a new class of \textit{mean field control game} (MFCG) problems arising naturally in the context of many large groups where agents are cooperating within each group but in competition with all agents in other groups. In this type of games, a MFC problem is defined at each group level motivating the dependency on the groups' distribution of the agents. At the full system level, a MFG problem is defined between groups explaining the freezing of the full system distribution and the following fixed point problem typical of this framework. Our algorithm naturally involves three learning rates: a fast one for the distribution of the group, a medium one for the agent's Q-matrix, and a slow one for the distribution of the overall population. We illustrate its performance on  linear-quadratic examples for which we derive explicit solutions for the optimal strategy.

In \cite{AFL2022}, the unified reinforcement Q-learning algorithm proposed in \cite{AFL2021} is generalized to finite horizon extended MFC and MFG problems. It is applied to the problem of a trader who wants to minimize transaction and inventory costs when trading an asset impacted by all agents' trades. We show in this paper that the algorithm can be naturally adapted for solving MFCG when both the distributions of  states and controls (for the group and for the overall population) are involved.

In Section \ref{sec:motiv}, we motivate the introduction of the new MFCG problem in the classical context of discrete time, finite horizon, differential games in discrete state and action spaces. Agents control their drifts and minimize an expected cost which may depend on the distributions of their own group and of the entire population. We give an intuitive justification of the fact that the solution of the MFCG provides an approximate Nash equilibrium between groups.

In Section \ref{sec:asymptotic} we introduce the discrete time and space infinite horizon MFCG considered in  this paper. We focus on the asymptotic formulation of the problem introduced in \cite{AFL2021} where comparisons with time-dependent and stationary formulations were discussed. In Section \ref{sec:extended}, we generalize the results of \cite{AFL2022} regarding finite time extended MFC and MFG problems to the MFCG framework. 

The  Q-learning approach to solve these problems is described in Section \ref{sec:Q-learning}
where we state the Bellman equation for the optimal action-value function (Q-matrix), and we introduce its three timescales stochastic approximation based on well-separated three learning rates: one for the states' distribution of the group, one for the action-value Q-matrix, and one for the the states' distribution of  the overall population.

Algorithm and learning rates are presented in Section \ref{sec:algo}. Its performance on a linear-quadratic benchmark are shown in Section \ref{sec:LQ-benchmark}. Section \ref{sec:traderPB} illustrates the results on the trader's problem, an example where the states' distribution of the group and the controls' distribution of the overall population appear in the objective function of the agent. We compare the strategies learned by our algorithm with the theoretical solutions provided in the Appendix \ref{sec:appendix_solutions}. 

\section{Mean Field Control Games}\label{sec:MFCG}

\subsection{Motivation} \label{sec:motiv}
In order to introduce our notion of MFCGs, we consider the familiar context of discrete time finite horizon differential games between agents evolving in a finite state space $\mathcal{X}$ by taking actions in a finite action space $\mathcal{A}$. The population is made of $M$ groups, each of size $N$. An agent will be indexed by a pair $(m,n)$ where the first index $m=1,\dots, M$ indicates her group and index $n=1,\dots, N$ being her identifier in the group. Agents are collaborating within their groups and competing with all  agents of  other groups. In other words,  all $N$ agents of group $m$ 
will try to collectively minimize the total cost of group $m$. Between groups, agents play a Nash equilibrium. So every single agent $(m,n)$ interacts, possibly in different ways, on both the distribution within its group and the distribution within the whole population.

We now present the general model that we consider, starting with the dynamics. At time $t=0,1,\dots,T-1$, agent $(m,n)$ uses the control $\alpha^{m,n}_{t} \in \mathcal{A}$. The evolution of her state is given by: $X^{m,n}_0 \sim \mu_0$ and for $t=0,1,\dots,T-1$,
\[
    \mathbb{P}(X^{m,n}_{t+1} = x' | X^{m,n}_{t} = x, \alpha^{m,n}_{t} = a, \mu_t = \mu) = p(x'|x,a,\mu),
\]
where $x$ and $x' \in \mathcal{X}$ represent respectively current and  next state, $a \in \mathcal{A}$ is the action taken, and $\mu \in \Delta^{|\X|}$ represents the empirical distribution of the whole population. Here, $p:\X \times \A \times \Delta^{ | \X |} \rightarrow \Delta^{ | \X |}$ is a transition kernel interpreted also as a function 
$$
    p:\X \times \X \times \A \times \Delta^{ | \X |} \rightarrow [0,1], \;\;\;\;\;\;\; (x,x',a,\mu) \mapsto p(x' | x,a,\mu),
$$ 
which provides the probability to jump to state $x'$ from state $x$ if action $a$ is taken. We assume that this transition kernel depends on the global distribution $\mu $ but not on the local distribution $\mu^m$ (that is the distribution of the agent's group defined below). 
In this finite-horizon setting, we allow for time-dependent feedback Markovian controls $\alpha : \{0,1,\dots,T-1\} \times \X \rightarrow \A$ that depend only on time and the state. So if agent $(m,n)$ uses control $\alpha$, then $\alpha^{m,n}_{t} = \alpha(t, X^{m,n}_{t})$. 

Considering the behavior of other groups as fixed, the goal for agents of group $m$ is to minimize the expected cost of the group:
\[
    J^m(\mathbf\alpha)=\frac{1}{N}\sum_{n=1}^N\mathds{E}\left[ \sum_{t=0}^T
    f(t,X_t^{m,n},\alpha_t^{m,n},\mu_t,\mu_t^m) +g(X_T^{m,n})\right],
\]
where $f$ is a running cost which may depend on the empirical distribution of the full population $\mu_t=\frac{1}{MN}\sum_{m=1}^M\sum_{n=1}^N \delta_{X^{m,n}_t}$
referred to as the \textit{global population}, and on the empirical distribution $\mu_t^m=\frac{1}{N}\sum_{n=1}^N \delta_{X^{m,n}_t}$ of the group $m$ 
referred to as the \textit{local population}. The terminal cost $g$ could as well depend on these empirical distributions at terminal time. Note that in this setting agents of group $m$ interact with agents of other groups through the distribution of the global population appearing in the cost.

Additional assumptions on $f$ and $g$ are needed, but we may keep in mind the simple quadratic cost case with, for example, $f(t,x,\alpha, \mu,\tilde\mu)= \frac{1}{2}\alpha^2 +\frac{c_1}{2}(x-\bar\mu)^2+ \frac{c_2}{2}\bar{\tilde\mu}^2$ where $\bar\mu$ and $\bar{\tilde\mu}$ denote respectively the means of the global population $\mu$ and the local population $\tilde\mu$. The first term is the classical quadratic cost for controlling the drift, the second term is an incentive to stay close to the global mean and the third term is a group incentive to keep the local mean close to zero. For simplicity we assume a zero terminal cost ($g=0$) in this example. We will revisit a linear-quadratic (LQ) continuous-time variant of this setting in Section~\ref{sec:LQ-benchmark}. The key point is that the interaction through the global mean is of mean field game (MFG) competitive nature, while the interaction through the local mean is of mean field control (MFC) collaborative nature, motivating the name mean field control game (MFCG). In other words, this problem is a competition between $M$ coalitions of $N$ players, all the players being identical in the sense that that they have similar dynamics and cost functions. The explicit solution for a continuous time and space version of this finite-player MFCG is given in Appendix~\ref{appendix}.

Passing to the mean field limit $M\to\infty, N\to\infty$ in a sense made precise in Appendix \ref{sec:approx}, a representative agent faces the following problem. Given a sequence of probability distributions $\boldsymbol{\mu}=(\mu_t)_{0\leq t \leq T}$, the goal is to solve the McKean-Vlasov (MKV) control problem of finding a minimizer $\hat\alpha$ for
\[
    J(\alpha)=\mathds{E}\left[ \sum_{t=0}^T
f(t,X_t^{\alpha,\boldsymbol{\mu}},\alpha_t,\mu_t,{\mathcal L}(X_t^{\alpha,\boldsymbol{\mu}})) +g(X_T^{\alpha,\boldsymbol{\mu}})\right],
\]
subject to:
\[
    X^{\alpha,\boldsymbol{\mu}}_{t+1} \sim p(X^{\alpha,\boldsymbol{\mu}}_{t}, \alpha_{t}, \mu_t), \quad t=0,1,\dots, T-1, \qquad X^{\alpha,\boldsymbol{\mu}}_0 \sim \mu_0.
\]
Allowing for time-dependent feedback Markovian controls means that $\alpha_t$ is given in the form $\alpha_t = \alpha(t, X^{\alpha,\boldsymbol{\mu}}_t)$ for some control function $\alpha$. Then, to find the Nash equilibrium, we need to solve the fixed point compatibility condition:
\[
    \mu_t = {\mathcal L} (X^{\hat\alpha,\boldsymbol{\mu}}_t),\quad \forall t\in \{0,1, \dots, T\},
\]
where $\mathcal{L}(X)$ denotes the law of the random variable $X$. This problem can be viewed as an MFG in which each player is of McKean-Vlasov type, in the sense that her dynamics and her cost function depend on her own distribution. As such, this can correspond to the limit of a MFTG~\cite{tembine} when the number of players goes to infinity. 
Solving this MFCG is justified by showing that the control $\hat\alpha$ enables the agents in the mixed finite-player game to achieve an $\epsilon$-Nash equilibrium. This argument is developed in the Appendix \ref{appendix} for the LQ example.

The proof of this result in a general setting  will be presented in the companion paper  \cite{MultiTimescaleConvergenceproof} in preparation. In particular, the analysis covers the case where the global distribution is involved in the dynamics. However, proving convergence when the local distribution appears in the dynamics is more challenging, which is why we do not include it in the dynamics studied in the present work.  The algorithm presented in this paper is in the context of finite state and action spaces. A version for continuous spaces based on a Deep Learning Actor Critic algorithm is a work in progress \cite{ADFLcontinuous}.

\subsection{Asymptotic Formulation}
\label{sec:asymptotic}
In this section we present the discrete time infinite horizon setting and we consider the asymptotic formulation of the game introduced in~\cite{AFL2021}. Our model involves the distribution of states within the collaborative agent's group (also called local distribution), and the distribution of states of the overall competitive population (also called global distribution).

We allow for time homogeneous controls $\alpha : \X \rightarrow \A$ that depend only on the state. We denote by $\mu^\alpha$ the asymptotic (long time) distribution of the controlled process following the strategy $\alpha$ which we assume to exist and to be unique (the state space being finite, aperiodicity and irreducibility of the discrete time process ensure these properties).

We go from finite horizon to infinite horizon so that the problem will be simpler to tackle with RL and we will look for stationary policy, see Section~\ref{sec:extended}. Given a cost function $f$ defined on $\X  \times \A \times \Delta^{ | \X |} \times \Delta^{ | \X |}$ and a discount rate $\gamma<1$,  we now consider the following infinite horizon asymptotic MFCG problem: 

Find a strategy $\hat{\alpha}$ and a distribution $\hat{\mu}$ such that:
\begin{enumerate}
\item (best response) $\hat{\alpha}$ is the minimizer of 
\begin{align*}
J(\alpha; \hat{\mu }) & =\mathds{E}\left[\sum_{t=0}^{\infty}\gamma^{t}f\left(X_{t}^{\alpha,\hat{\mu}},\alpha(X_{t}^{\alpha,\hat{\mu} }),\hat{\mu},\mu^{\alpha, \hat\mu}\right)\right],
\end{align*}
where $X_{0}^{\alpha,\hat\mu}\sim\mu_{0}$ and for $t=0,1,\dots$,
\begin{align*}
    \mathbb{P}(X_{t+1}^{\alpha,\hat{\mu}} = x' | X_{t}^{\alpha,\hat{\mu}} = x, \alpha(X_{t}^{\alpha,\hat{\mu}}) = a, {\mu} = \hat\mu) = p(x'|x,a,\hat\mu).
\end{align*}
and $\mu^{\alpha, \hat\mu}=\lim_{t\rightarrow\infty}\mathcal{L} (X_{t}^{\alpha,\hat{\mu}} )$.
\item (fixed-point) $\hat{\mu}=\lim_{t\rightarrow\infty}\mathcal{L}(X_{t}^{\hat\alpha,\hat{\mu}})=\mu^{\hat{\alpha}, \hat\mu}$.
\end{enumerate}

In order to make sense of the above problem statement we have to restrict to actions $\alpha : \X \rightarrow \A$ which are such that the controlled process $X_{t}^{\alpha,\mu}$ has a limiting distribution, {\textit i.e.}, $\lim_{t\rightarrow\infty}\mathcal{L} (X_{t}^{\alpha, \mu} )$ exists. For a finite state Markov chain this is the case if $( X_{t}^{\alpha, \mu})_{t\in \mathbb{N}}$ is irreducible and aperiodic. We therefore assume that the strategy $\hat{\alpha}$ is the minimizer over all strategies such that $( X_{t}^{\alpha, \hat{\mu}})_{t\in \mathbb{N}}$ is irreducible and aperiodic.

\remark{\label{rem:statio-asympto}
We could have also considered a classical formulation of MFCG where $\hat\mu$ and $\mu^{\alpha, \hat\mu}$ are flows of distributions $(\hat\mu_t)_{t\in \mathbb{N}}$ and $(\mu_t^{\alpha, \hat\mu})_{t\in \mathbb{N}}$ in which case the fixed-point requirement is $\hat\mu_t=\mathcal{L}(X_{t}^{\alpha,\hat{\mu}})$ for every $t\in \mathbb{N}$, and where the strategy is time-dependent. As well, we could have considered a stationary formulation of MFCG where $\mu^\alpha$ is the stationary distribution of the controlled process, equal to $\hat\mu$ in the fixed-point step. As in \cite{AFL2021}, it can be shown that the optimal strategies for the asymptotic and stationary problems coincide, and they coincide with the limiting optimal strategy (as $t\to\infty$) of the classical formulation.
}

\subsection{Finite Horizon Extended Formulation}\label{sec:extended}
Following \cite{AFL2022}, we generalize the MFCG problem and its reinforcement learning algorithm to the case with a discrete time finite horizon $T$, on a finite state space, and mean field of state and control. The state-action space is as described in Section \ref{sec:MFCG}. The state follows a random evolution in which $X_{t+1}$ is determined as a function of the current state $X_t$, the action $\alpha_t$, and some noise. We introduce the transition probability function:
$$
    p(x'|x,a,\nu), \qquad (x, x', a, \nu) \in \mathcal{X} \times \mathcal{X} \times \mathcal{A} \times \Delta^{|\mathcal{X} \times \mathcal{A}|}, 
$$
which provides the probability to jump to state $x'$ given its current state $x$, the action taken $a$ and the global population distributed as $\nu$. We assume no dependence on the state-action group distribution $\tilde\nu$ in order to apply the MKV Bellman equation introduced in \cite{AFL2021}. For simplicity, we consider the homogeneous case where this function does not depend on time.
Restoring this time-dependence if needed is a straightforward procedure.

We now consider the MFCG cost function given by: for $\nu = (\nu_t)_{t=0,1,\dots,T}$ 
$$
    J(\alpha;\nu) = \mathbb{E}\left[ \sum_{t=0}^{T-1} f(X^{\alpha,\nu}_t, \alpha_t, \nu_t, \nu^{\alpha, \nu}_t) + g(X^{\alpha,\nu}_{T}, \mu_{T}, \mu^{\alpha, \nu}_{T}) \right],
$$
where 
$\mu_{T}$ (resp. $\mu^{\alpha, \nu}_{T}$) is the first marginal of $\nu_{T}$ (resp. $\nu^{\alpha, \nu}_{T})$. Again, for simplicity, we assume that $f$ does not depend on time. The process $X^{\alpha,\nu}$ has a given initial distribution $\mu_0 \in \Delta^{|\mathcal{X}|}$ and follows the dynamics 
$$
    \mathbb{P}(X^{\alpha,\nu}_{t+1} = x' | X^{\alpha,\nu}_t = x, \alpha_t = a, \nu_t = \nu) = p(x'|x,a, \nu).
$$

\section{Q-Learning}\label{sec:Q-learning}

\subsection{Action-Value Function}

Our algorithm to solve the MFCG is based on the concept of $Q$-learning which is a well known procedure to solve Markov decision problems. However, following \cite{AFL2021} we combine the idea of $Q$-learning with the model agnostic view of reinforcement learning. We first adapt the $Q$-learning concepts to our problem at hand. Since the local distribution is not fixed and depends on the control itself, we have to adapt the classical $Q$-learning in the spirit of \cite{AFL2021}. For an admissible control $\alpha : \X \rightarrow \A$ and a pair $(x,a) \in \X \times \A$, we define the new control $\alpha_{x,a}$ by 
\begin{equation}
\alpha_{x,a} (x')=\begin{cases}
      a & \text{if   } x' = x, \\
      \alpha(x')  & \text{otherwise.}
\end{cases}\label{mod:alpha}
\end{equation}
Given a global measure $\mu$ and a strategy $\alpha$, the $Q$-function for our problem is given by:
\[
Q^{\alpha}_{\mu} (x,a)=f(x,a,\mu,\mu^{\alpha_{x,a}, \nu}) + \mathbb{E}\left[\sum_{t=1}^{\infty}\gamma^{t}f(X_{t},\alpha(X^{\alpha,\mu}_{t}),\mu,\mu^{\alpha, \nu})\lvert X^{\alpha,\mu}_{0}=x,A_{0}=a\right],
\]
where $\mu^{\alpha_{x,a}, \nu}$ is the local distribution relative to the strategy $\alpha_{x,a}$.
The optimal function in the sense of minimizing cost is given by
$$Q_{\mu}^{*} (x,a) := \min_{\alpha} Q_{\mu}^{\alpha}(x,a).$$ 
From the function $Q_{\mu}^{*}$ one obtains the control $\alpha^* (x)=\arg \min_a Q_{\mu}^{*} (x,a)$ (in fact in the algorithm presented in Section \ref{sec:algo}, we use a randomized policy, which is not taken into account here). Note that the minimizing strategy may depend on the global measure $\mu$. For fixed $\mu$, the function $Q_{\mu}^{*}$ follows the Bellman equation given by: 
\begin{equation}\label{MKV-Bellman}
Q_{\mu}^{*}(x,a)=f(x,a,\mu,\mu_{x,a}^{*,\mu})+\gamma\sum_{x'}p(x'\lvert x,a,\mu)\min_{a'}Q_{\mu}^{*}(x',a').
\end{equation}
Note that using this modified (McKean--Vlasov type) Bellman equation established in \cite{AFL2021} allows us to consider the Q-function as a function of state and action only.
The measure $\mu_{x,a}^{*, \mu} = \lim_{t\rightarrow\infty}\mathcal{L} (X_{t}^{\alpha^*_{x,a},\mu})$ corresponds to the strategy $\alpha^*_{x,a}$ which is derived from $\alpha^*$ by changing the action in state $x$ to $a$, see \eqref{mod:alpha}.
The above Bellman equation follows from the results in \cite{AFL2022} as the measure $\mu$ is fixed and does not depend on $\alpha$.

\subsection{Time-Dependent Q-Function}
The definition of the time-dependent optimal $Q$-function in the extended framework is given for a fixed flow of state-action global distributions $\nu=(\nu_t)_{t=0,1,\dots,T}$ by:
 \begin{equation*}%
 \left\{
 \begin{aligned}
 &Q^*_{T,\nu}(x,a) = g(x,\mu_{T}, \mu^{\alpha, \nu}_{T}), \qquad (x,a) \in \mathcal{X} \times \mathcal{A},
 \\
    &Q^*_{t,\nu}(x,a)=\min_\alpha\mathbb{E}\left[\sum_{t'=t}^{T-1}  f(X^{\alpha, \nu}_{t'},\alpha_{t'}(X^{\alpha, \nu}_{t'}),\nu_{t'}, \nu^{\alpha, \nu}_{t'}) + g(X^{\alpha, \nu}_{T},\mu_{T}, \mu^{\alpha, \nu}_{T})\,\Big\vert\, X^{\alpha, \nu}_t=x,A_t=a\right], 
    \\
    &\qquad\qquad t=0,1,\dots,T-1, \quad (x,a) \in \mathcal{X} \times \mathcal{A} ,
    \end{aligned}
    \right.
\end{equation*}
where $\mu_{T}$ (resp. $\mu^{\alpha, \nu}_{T}$) is the first marginal of $\nu_{T}$ (resp. $\nu^{\alpha, \nu}_{T})$, and $\alpha_{t'}(\cdot) = \alpha({t'},\cdot)$. Using dynamic programming, it can be shown that $(Q^*_{t,\nu})_{t}$ is the solution of the Bellman equation:  
 \begin{equation*}%
 \left\{
 \begin{aligned}
 &Q^*_{T}(x,a) = g(x,\mu_{T}, \mu^{\alpha, \nu}_{T}), \qquad (x,a) \in \mathcal{X} \times \mathcal{A},
    \\
   & Q^*_{t,\nu}(x,a) = f(x, a,\nu_t,\nu_t^{\tilde\alpha, \nu}) + \sum_{x' \in \mathcal{X}} p(x' | x, a,\nu) \min_{a'} Q^*_{t+1,\nu}(x',a'), \qquad t=0,1,\dots,T-1,\quad (x,a) \in \mathcal{X} \times \mathcal{A},
    \end{aligned}
    \right.
\end{equation*}
where $\nu_t^{\tilde\alpha, \nu}$ takes into account the modification of $\alpha$ due to the decision $a$ at state $x$.
The corresponding optimal value function $(V_{t,\nu}^*)$ is given by:
$$
    V^*_{t,\nu}(x) = \min_a Q^*_{t,\nu}(x,a), \qquad t=0,1,\dots,T,\quad x \in \mathcal{X}.
$$
One of the main advantages of computing the action-value function instead of the value function is that from the former one obtains the optimal control at time $t$ by computing $\arg\min_{a \in \mathcal{A}} Q^*_t(x,a)$. This is particularly important in order to design model-free methods as we will see in the next section.

The next step consists in describing the updates of the $Q_t$'s tables, the flows of measures $\nu_t$'s and $\nu^{\alpha, \nu}_t$'s. As for the infinite horizon case discussed in the next section, a three-timescale approach is implemented by introducing three learning rates $\rho_k^\nu<\rho_k^{Q}<\rho_k^{\nu^{\alpha, \nu}}$. We skip the details for the finite horizon extended framework as they are similar to \cite{AFL2022} where it is presented for the two timescale case. This approach justifies the algorithm presented in Section \ref{sec:algo}.

\subsection{Stochastic Approximation}\label{sec:rates}

In this section we propose a learning procedure that under reasonable assumptions on the functions $p$ and $f$ approximates the solution of the discrete time MFCG. The algorithm is based on the idea that the local distribution, the $Q$-function describing the optimal strategy, and the global distribution should be updated at different rates. For the sake of a lighter notation, we will use the notation $\mu, Q$ and $\mu^\alpha$ omitting the mutual dependencies that are fully discussed in the previous sections.

For a pure MFC and a pure MFG problem  the authors of \cite{AFL2021} use results in \cite{borkar1997stochastic,MR2442439} for classical $Q$-learning to show that a two-timescale approach involving the system distribution and the optimal response can converge to either the MFC solution or the MFG solution depending on how the learning rates are chosen. For a MFC problem, the system distribution resulting from a chosen strategy has to be updated more frequently than the strategy itself. In contrast, the MFG case requires the strategy to be updated more frequently than the distribution. 

To gain some intuition for the three-timescale approach used to approximate our MFCG, we start with the function $Q:\X\times \A \rightarrow \mathbb{R}$ inducing a strategy $\alpha'$ such that at each time the system is at state $x$, the action $\arg \min_a Q (x,a)$ is chosen. Say the global distribution $\mu$ is frozen (as it is in part 1. of our MFCG problem) and the local distribution is given by $\mu^{\alpha}$, then the local population will be driven at the next step towards the new distribution $\sum_{x\in \X}\mu^{\alpha}(x) p(x'\lvert x,\arg\min_{a}Q(x,a),\mu)$ if all players follow the strategy encoded in $Q$. This continues until a fixed-point $\mu^{\alpha'}$ is reached. When a fixed-point is (approximately) reached, the strategy has to be updated, taking this new limiting distribution into account. This leads to a new optimal strategy with action values given by
\[
f(x,a,\mu,\mu^{\alpha'})+\gamma\sum_{x'}p(x'\lvert x,a, \mu)\min_{a'}Q(x',a').
\]
This procedure continues until an optimal pair of strategy $\alpha$ and resulting limiting measure $\mu^{\alpha}$ is reached which depends on the frozen global measure $\mu$. In an outer global optimization the fixed-point for the global measure is now obtained by updating the global measure via $\sum_{x\in \X}\mu(x) p(x'\lvert x,\arg\min_{a}Q(x,a),\mu)$. The three timescales therefore arise naturally by the different layers of optimization involved in the problem. It is intuitive that in each layer one has to perform sufficiently many iterations to ensure that the optimization in the next layer is based on sufficiently accurate results. This idea leads to a learning rate that decreases from the outer to the inner layer. In addition the ratios of the increasing learning rates (from inner to outer layer) have to be sufficiently large.  

These considerations lead to the following updating rules:
$(\mu,Q,\mu^{\alpha})$ are updated with rates $\rho_{k}^{\mu}<\rho_{k}^{Q}<\rho_{k}^{\mu^{\alpha}}$by 
\begin{equation}\label{update:system}
\begin{cases}
\mu_{k+1}=\mu_{k}+\rho_{k}^{\mu}\mathcal{P}(\mu_{k},Q_{k},\mu_{k}),\\
Q_{k+1}=Q_{k}+\rho_{k}^{Q}\mathcal{T}(\mu_{k},Q_{k},\mu^{\alpha}_{k}),\\
\mu^{\alpha}_{k+1}=\mu^{\alpha}_{k}+\rho_{k}^{\mu^{\alpha}}\mathcal{P}(\mu_{k},Q_{k},\mu^{\alpha}_{k}),
\end{cases}
\end{equation}
where $k$ denotes the learning episode (see the algorithm below), and
\[
\begin{cases}
\mathcal{P}(\mu,Q,\nu)(x)=(\nu P^{\mu,Q})(x)-\nu(x),\\
\mathcal{T}(\mu,Q,\mu^{\alpha})(x,a)=f(x,a,\mu,\mu^{\alpha})+\gamma\sum_{x'}p(x'\lvert x,a,\mu)\min_{a'}Q(x',a')-Q(x,a)\\
P^{\mu,Q}(x,x')=p(x'\lvert x,\arg\min_{a}Q(x,a),\mu),\\
(\nu P^{\mu,Q})(x)=\sum_{x_{0}}\nu(x_{0})P^{\mu,Q}(x_{0},x).
\end{cases}
\]
To see that the above system does in fact converge to a solution of our MFCG, 
we assume that $\rho^{\mu}_k\ll \rho^{Q}_k$ so that
$\rho^{\mu}_k / \rho^{Q}_k$ is of order $\epsilon \ll 1$, and $\rho^{Q}_k\ll \rho^{\mu^{\alpha}}_k$ so that 
$\rho^{Q}_k/\rho^{\mu^{\alpha}}_k$ is of order $\tilde\epsilon \ll 1$.

Now, following \cite{borkar1997stochastic}, we denote by $\tau$ a continuous time variable, and we consider the following ODE system which tracks the system \eqref{update:system}:
\[
\begin{cases}\label{ODE:system}
\dot{\mu}_\tau = \mathcal{P}(\mu_{\tau},Q_{\tau},\mu_{\tau})\\
\dot{Q}_{\tau}=\frac{1}{\epsilon}\mathcal{T}(\mu_{\tau},Q_{\tau},\mu^{\alpha}_{\tau}) \\
\dot{\mu}^{\alpha}_{\tau}=\frac{1}{\epsilon \cdot \tilde{\epsilon}} \mathcal{P}(\mu_{\tau},Q_{\tau},\mu^{\alpha}_{\tau}).
\end{cases}
\]
Furthermore, we assume that the functions $f$ and $p$ are such that the system fulfills a Lipschitz condition. As shown in \cite{AFL2021}, this can be ensured by Lipschitz continuity of $f$ and $p$ and by smoothing the minimum in the definition of $P$. We refer to \cite{AFL2021} where these considerations are treated in more detail for a two-timescale approach. 

We start with the fastest timescale. For a fixed global distribution $\mu$ and a fixed action table $Q$, we assume that the ODE 
\[ 
\dot{\mu}^{\alpha}_{\tau}=\frac{1}{\epsilon \cdot \tilde{\epsilon}} \mathcal{P}(\mu,Q,\mu^{\alpha}_{\tau})
\]
has a unique asymptotically {($\tilde\epsilon\rightarrow 0$)} stable equilibrium $\mu^{Q,\mu}$ such that $\mathcal{P}(\mu, Q,\mu^{Q,\mu})=0$. Now, we plug this equilibrium $\mu^{Q, \mu}$ into the second equation and obtain the ODE 
\[ 
\dot{Q}_{\tau}=\frac{1}{\epsilon}\mathcal{T}(\mu,Q_{\tau},\mu^{Q_\tau, \mu}).
\]
Again, we assume that the above ODE has a stable equilibrium ($\epsilon\to 0$), which we call $Q^\mu$ and which satisfies that $\mathcal{T}(\mu,Q^\mu,\mu^{{Q^\mu},\mu})=0$. 
Now, going to the slowest timescale, the first equation has an asymptotic ($\tau\to\infty)$ equilibrium, say $\mu_\infty$ that solves  $\mathcal{P}(\mu_{\infty},Q^{\mu_\infty},\mu_{\infty})=0$. By uniqueness of this equilibrium, we get that $\mu_\infty=\mu^{ Q^{\mu_\infty}, \mu^{\infty}}$, 
which in turns implies that $\mu_{\infty}$ and the action given by minimizing $Q^{\mu_\infty}$ solves our MFCG.

\section{Reinforcement Learning Algorithm}\label{sec:algo}

\subsection{Asymptotic Version}
The three-timescale mean field Q-learning algorithm (U3-MF-QL) that we propose leverages the two-timescale version (U2-MF-QL) introduced by \cite{AFL2021}. It not only encompasses learning the pure MFG and pure MFC problems, but, more importantly, it facilitates learning the generalized MFCG problems. Despite of its advantage and flexibility, it inherits the very simple intuition that by manipulating the relative value of learning rates we can induce the algorithm to updating distributions in either MFG's or MFC's manner, as described in the last Section~\ref{sec:rates}. Depending on whether the problem has infinite horizon or finite horizon, the U3-MF-QL algorithm will be specified accordingly. Here we first introduce the infinite horizon version (U3-MF-QL-IH) in Algorithm~\ref{alg:rl_infinite}. The intuition underlying the algorithm is based on the asymptotic formulation but, as explained in Remark~\ref{rem:statio-asympto}, this is equivalent to solving the stationary problem. Furthermore, for the sake of simplicity we present the algorithm for an MFCG involving only the state distribution but it can be adapted to solve an MFCG involving the state-action distribution, i.e., an extended MFCG. In Section \ref{sec:algo_finite}, Algorithm~\ref{alg:rl_finite} is presented for the finite horizon extended MFCG. Both algorithms can be adapted to solve continuous states and actions problems by applying the necessary truncation and discretization techniques as originally discussed in \cite{AFL2021}.

\begin{algorithm}[H]
\caption{Three-Timescale Mean Field Q-Learning - Discrete Time Infinite Horizon}
\begin{algorithmic}[1]
\Require{\\
    \qquad Time steps $t=0,1,\dots,T$ with $T>>0$, \\
    \qquad Finite state space: $\X = \{x_0, \dots, x_{|\X|-1}\}$, \\
    \qquad Finite action space: $\A = \{a_0, \dots, a_{|\A|-1}\}$, \\
    \qquad Initial distribution of the representative player: $\mu_0$, \\
    \qquad Factor of the $\varepsilon$-greedy policy: $\varepsilon$, \\
    \qquad Break rule tolerances: $tol_Q$, $tol_{\mu}$, $tol_{\tilde{\mu}}$.
}
\State{\textbf{Initialization:}}
    \State{\qquad $Q^0(x, a)=0$ for all $(x,a) \in \X \times \A$,}
    \State{\qquad $\mu^0_{t} = \frac{1}{|{\X}|}J_{|{\X}|}$ and $\tilde{\mu}^0_{t} = \frac{1}{|{\X}|}J_{|{\X}|}$ for $t \le T$,}
    \State{\qquad where $J_{d}$ is a $d$-dimensional vector.}
\For{each episode $k=1, 2, \dots$}
    \State{\textbf{Observe initial state:} $X_0^k \sim \mu_T^{k-1}$ and set $Q^{k} \equiv Q^{k-1}$.}
    \For{$t = 0, \dots, T$}
        \State{\textbf{Choose action:}}
            \State{\qquad choose $A_t^k$ using the $\varepsilon$-greedy policy derived from $Q^{k}(X_t^k, \cdot)$.}
        \State{\textbf{Update distributions:}}
            \State{\qquad $\mu^{k}_t = \mu^{k-1}_t + \rho_{k}^{\mu}(\boldsymbol{\delta}(X_t^k) - \mu^{k-1}_t)$,}
            \State{\qquad $\tilde{\mu}^{k}_t = \tilde{\mu}^{k-1}_t + \rho_{k}^{\tilde{\mu}}(\boldsymbol{\delta}(X_t^k) - \tilde{\mu}^{k-1}_t)$,}
            \State{\qquad where $\boldsymbol{\delta}(X_t^k)=\left(\mathbf{1}_{x}(X_t^k)\right)_{x \in \X}$.}
        \State{\textbf{Observe next state:}}
            \State{\qquad observe $X_{t+1}^k$ from the environment.}
        \State{\textbf{Observe cost:}}
            \State{\qquad observe $f_t = f(X_t^k, A_t^k, \mu^{k}_t, \tilde{\mu}^{k}_t)$.}
        \State{\textbf{Update $Q$ table:}}
            \State{\qquad $Q^{k}(x,a) = Q^{k}(x,a) +  \mathbf{1}_{x,a}(X_t^k, A_t^k)\rho^{Q}_{x,a, t, k}\left(f_t + \gamma \min_{a' \in \A} Q^{k}(X_{t+1}^k, a')-Q^{k}(x,a)\right)$,}
         \State{\qquad where $\gamma$ is the discount parameter.}
    \EndFor
    \If{$\norm{Q^{k}-Q^{k-1}}\le tol_{Q}$, $\norm{\mu^{k}-\mu^{k-1}}\le tol_{\mu}$, and $\norm{\tilde{\mu}^{k}-\tilde{\mu}^{k-1}}\le tol_{\tilde{\mu}}$ }
        \State{break}
    \EndIf
\EndFor
\end{algorithmic}
\label{alg:rl_infinite}
\end{algorithm}

\subsubsection{Learning Rates}
By choosing $\rho^{\mu}_k < \rho^{Q}_k$, we induce the global distribution $\mu$ to converge in the fashion of MFG. On the other hand, by letting $\rho^{Q}_k < \rho^{\tilde{\mu}}_k$, we allow the local distribution $\tilde{\mu}$ to renew towards the MFC style. Combining both such that $\rho^{\mu}_k < \rho^{Q}_k < \rho^{\tilde{\mu}}_k$, the algorithm is expected to learn both the global and local distributions simultaneously. In addition, to ensure that the learned Q-table and distributions can stablize at the end of the episode iteration, all the three learning rates shall also decay as the number of episodes, $k$, increases. Adapting the learning rate discussed in \cite{AFL2021}, we design the triplet of learning rates as follows.
\begin{align}\label{eq:rates}
\rho_{x,a, t, k}^{Q}:=\frac{1}{(1+ \#\abs{(x, a, k, t)})^{\omega^Q}}, &&
\rho_{k}^{\mu}:=\frac{1}{(1+k)^{\omega^{\mu}}}, &&
\rho_{k}^{\tilde{\mu}}:=\frac{1}{(1+k)^{\omega^{\tilde{\mu}}}},
\end{align}
where $\#\abs{(x, a, k, t)}$ counts the visits of the pair $(x,a)$ up to the episode $k$ and time $t$. The triplet $(\omega^{Q}, \omega^{\mu}, \omega^{\tilde{\mu}})$ should be chosen such that $\omega^{\mu} > \omega^{Q} >\omega^{\tilde{\mu}}$, so that $\rho^{\mu}_k < \rho^Q_k < \rho^{\tilde{\mu}}_k$, and it should satisfy $\omega^{Q} \in (0.5, 1)$. 

\subsection{Time-Dependent Version}\label{sec:algo_finite}
The three-timescale mean field Q-learning approach specified for the finite horizon (U3-MF-QL-FH) is shown in Algorithm~\ref{alg:rl_finite}. Although its overall structure is similar to that of Algorithm~\ref{alg:rl_infinite}, we shall highlight several important differences. First, in the finite horizon problem, the algorithm must learn the optimal control and state-action distribution for each time point. So, the number of $Q$ tables to be learned is $T -1$, each corresponding to a time step, except for the terminal time which is excluded because no action is taken at time $T$. In contrast, in the infinite horizon problem we had just a single $Q$ table to learn. Second, in each episode, the initial state $X_0$ is always drawn from the initial distribution $\mu_0$. This is in contrast to the infinite horizon case where the initial state $X_0$ is drawn from the terminal empirical distribution learned up to the last episode $\mu^{k-1}_{T}$. Third, within each episode, the algorithm only iterates through the time steps from $0$ to ${T - 1}$. It skips the terminal time $T$, because at time ${T-1}$, once the action $A_{T-1}$ is chosen, the final state $X_{T}$ can be generated and henceforth the terminal cost $g(X_{T})$ is observed, which already completes the episode. In the infinite horizon case, whether one iterates up to ${{T-1}}$ or ${{T}}$ does not make a big difference. Fourth, when updating the $Q_t$ table, the table $Q_{t+1}$ for the next time step $t < {T}$ needs to be taken into account, in contrast to the infinite horizon case. Lastly, the learning rate for the $Q$-tables in the finite horizon case are 
\begin{align}
\rho_{x,a,k}^{Q_t}:=\frac{1}{(1+ T\#\abs{(x, a, k, t)})^{\omega^Q}},
\end{align}
where $\#\abs{(x, a, k, t)}$ counts separately for each time step $t$ the visits of tuples $(x,a)$ up to episode $k$. 

The approximation of this time-dependent version of the algorithm to the MFCG solution can be shown similarly as we did for the asymptotic problem in Section~\ref{sec:rates}. We refer the reader to \cite{AFL2021} where this is done for the pure mean field control and the pure mean field game problem. 

\begin{algorithm}[H]
\caption{Three-Timescale Mean Field Q-Learning - Discrete Time Finite Horizon}
\begin{algorithmic}[1]
\Require{\\
    \qquad Time steps: $t = 0,1, \dots, T$, \\
    \qquad Finite state space: $\X = \{x_0, \dots, x_{|{\X}|-1}\}$, \\
    \qquad Finite action space: $\A = \{a_0, \dots, a_{|{\A}|-1}\}$, \\
    \qquad Initial distribution of the representative player: $\mu_0$, \\
    \qquad Factor of the $\varepsilon$-greedy policy: $\varepsilon$,\\
    \qquad Break rule tolerances: $tol_Q$, $tol_{\nu}$, $tol_{\tilde{\nu}}$.
}
\State{\textbf{Initialization}:}
    \State{\qquad $Q^0_t(x, a)=0$ for all $(x,a) \in \X \times \A$, for all $t \in \{0,1,\dots,T\}$,}
    \State{\qquad $\nu^0_t = \frac{1}{\abs{\X \times \A}}J_{|{\X}| \times |{\A}|}$, $\tilde{\nu}^0_t = \frac{1}{\abs{\X \times \A}}J_{|{\X}| \times |{\A}|}$ for all $t \in \{0,1,\dots,T\}$,}
    \State{\qquad where $J_{n \times m}$ is an $n \times m$ unit matrix.}
\For{each episode $k=1, 2, \dots$}
    \State{\textbf{Observe initial state:} $X_0 \sim \mu_0$.}
    \For{$t = 0$,1, \dots, ${T-1}$}
        \State{\textbf{Choose action:}}
            \State{\qquad choose $A_t$ using the $\epsilon$-greedy policy derived from $Q^{k-1}_t(X_t, \cdot)$.}
        \State{\textbf{Update empirical distributions:}}
            \State{\qquad $\nu^{k}_t = \nu^{k-1}_t + \rho_{k}^{\nu}(\boldsymbol{\delta}(X_t, A_t) - \nu^{k-1}_t)$,}
            \State{\qquad $\tilde{\nu}^{k}_t = \tilde{\nu}^{k-1}_t + \rho_{k}^{\tilde{\nu}}(\boldsymbol{\delta}(X_t, A_t) - \tilde{\nu}^{k-1}_t)$,}
            \State{\qquad where $\boldsymbol{\delta}(X_t, A_t)=(\mathbf{1}_{x,a}(X_t, A_t))_{x \in \X, a \in \A}$.}
        \State{\textbf{Observe next state:}}
            \State{\qquad observe $X_{t+1}$ from the environment}
        \State{\textbf{Observe cost:}}
            \State{\qquad running cost $f_t = f(X_t, A_t, \nu^{k}_t, \tilde{\nu}^{k}_t)$,}
            \State{\qquad terminal cost $g_T = g(X_T)$ when reach $t+1=T$.}
        \State{\textbf{Update} $Q_t$:}
            \State{\qquad $Q^{k}_t(x,a) = Q^{k-1}_t(x,a) +  \mathbf{1}_{x,a}(X_t, A_t)\rho^{Q}_{x,a, k, t}\left(f_t + B-Q^{k-1}_t(x,a)\right)$,}
            \State{\qquad where $B = \mathbf{1}_{\{t+1=T\}} g_T + \mathbf{1}_{\{t+1 < T\}}\min_{a \in \A} Q^{k-1}_{t+1}(X_{t+1},a)$.}
    \EndFor
    \If{ $\norm{Q^{k}-Q^{k-1}}\le tol_{Q}$, $\norm{\nu^{k}-\nu^{k-1}}\le tol_{\nu}$, and $\norm{\tilde{\nu}^{k}-\tilde{\nu}^{k-1}}\le tol_{\tilde{\nu}}$ }
        \State{break}
    \EndIf
\EndFor
\end{algorithmic}
\label{alg:rl_finite}
\end{algorithm}

\section{Numerical Experiments}\label{sec: Numerical results}

We illustrate the performance of our algorithms on benchmark models for which we have explicit solutions:  in the infinite horizon case (Algorithm \ref{alg:rl_infinite}) in Section \ref{sec:LQ-benchmark}, and  in a finite horizon extended game setting (Algorithm \ref{alg:rl_finite}) in Section \ref{sec:traderPB}.

\subsection{Asymptotic Problem}\label{sec:LQ-benchmark}

For MFG or MFC problems in finite spaces, explicit solutions are usually not available and we would have to rely on approximate solutions obtained by other numerical methods to compare with the solutions obtained with our RL algorithm. Instead, here we choose to work with a linear-quadratic model in continuous time and space for which we can easily derive explicit solutions (see Appendix \ref{appendix:LQ}). We then apply our algorithm to a discretization in time and space of this model described in Section \ref{sec:results}. We do not address here the quality of this discretization which has been widely studied. We simply compare the results of our algorithm with the explicit solutions of the continuous model.

Specifically, we consider a continuous-time and space benchmark linear-quadratic MFCG problem with a running cost given by 
\begin{align}\label{eq:f}
f(x,\alpha,\mu,\mu^{\alpha, \mu}) & =\frac{1}{2}\alpha^{2}+c_{1}(x-c_{2} m)^{2}+c_{3}(x-c_{4})^{2}+\tilde{c}_{1}(x-\tilde{c}_{2}{m}^{\alpha, \mu})^{2}+\tilde{c}_{5}({m^{\alpha, \mu}})^{2},
\end{align}
where $m=\int x\mathrm{d}\mu(x)$ and ${m}^{\alpha, \mu}=\int x\mathrm{d}\mu^{\alpha, \mu}(x)$, and
$c_{1}$, $\tilde{c}_{1}$, and $\tilde{c}_{5}$ are positive
constants. Here, $\mu$ and $\mu^{\alpha, \mu}$ are understood as \textit{global
environment} and \textit{local environment}. The constant $c_{1}$ determines the magnitude of the \textit{global
effect} and the constants $\tilde{c}_{1}$, $\tilde{c}_{5}$ specify \textit{local effects}.

The asymptotic formulation of this MFCG problem is given
by
\begin{align*}
\inf_{\alpha}J(\alpha; \mu) & =\inf_{\alpha}\mathds{E}\left[\int_{0}^{\infty}\mathrm{e}^{-\beta t}f(X_{t}^{\alpha, \mu},\alpha_{t},\mu,\mu^{\alpha, \mu})\mathrm{d}t\right]\\
 & =\inf_{\alpha}\mathds{E}\left[\int_{0}^{\infty}\mathrm{e}^{-\beta t}\left(\frac{1}{2}\alpha_{t}^{2}+c_{1}(X_{t}^{\alpha, \mu}-c_{2}m)^{2}+c_{3}(X_{t}^{\alpha, \mu}-c_{4})^{2}+\tilde{c}_{1}(X_{t}^{\alpha, \mu}-\tilde{c}_{2}m^{\alpha, \mu})^{2}+\tilde{c}_{5}\left(m^{\alpha, \mu}\right)^{2}\right)\mathrm{d}t\right].
 \end{align*}
subject to $\mathrm{d}X_{t}^{\alpha, \mu}=\alpha_{t}\mathrm{d}t+\sigma\mathrm{d}W_{t},\, X_{0}^{\alpha, \mu}\sim\mu_{0}$,
and the fixed point condition
    $m = \lim_{t \rightarrow \infty} \mathds{E}(X^{{\hat{\alpha}, \mu}}_t)=m^{\hat\alpha,\mu}$,
where $\hat{\alpha}$ is the optimal action.

\subsubsection{Results}\label{sec:results}

We consider the asymptotic MFCG with the following choice of parameters: $c_1=0.5$, $c_2=1.5$, $c_3=0.5$, $c_4=0.25$ $\tilde{c}_1=0.3$, $\tilde{c}_2=1.25$, $\tilde{c}_5=0.25$, discount rate $\beta = 1$, and volatility of the state dynamics $\sigma = 0.5$. We truncate the infinite time horizon at $T=20$, and discretize the interval $[0,T]$ with time steps of size $\Delta t = 10^{-2}$. The discount factor in the discrete time setting is then given by $\gamma := e^{-\beta \Delta t}$. The state space is $\X = \{x_0=-2 + x_c, \dots, x_{|{\X}|-1}=2+x_c\}$ centered at $x_c=0.25$, and the action space is $\A = \{a_0 = -3, \dots, a_{|{\A}|}=3\}$, where the step sizes are $\Delta x = \Delta a = \sqrt{\Delta t}=10^{-1}$. The $\epsilon$-greedy policy parameter is $0.01$. We remind the reader on the choice of the learning rates. In contrast to the pure MFG and MFC problems which can be learned by the two-timescale parameterization proposed in \cite{AFL2021}, the MFCG problem requires the three-timescale as explained in Section~\ref{sec:rates}. We choose the three learning rates to be $(\omega^{\mu}, \omega^{Q}, \omega^{\mu^{\alpha}}) = (0.85, 0.55, 0.15)$, which satisfy $\rho_{k}^{\mu}<\rho_{k}^{Q}<\rho_{k}^{\mu^{\alpha}}$. In addition, we also demonstrate that if we miss-specify the learning rates $(\omega^{\mu}, \omega^{Q}, \omega^{\mu^{\alpha}})$ by either $(0.85, 0.55, 0.85)$ in which case $\rho_{k}^{\mu}=\rho_{k}^{\mu^{\alpha}}<\rho_{k}^{Q}$ or by $(0.15, 0.55, 0.15)$ in which case $\rho_{k}^{\mu}=\rho_{k}^{\mu^{\alpha}}>\rho_{k}^{Q}$, then the algorithm fails to learn the correct result.

\begin{figure}[H]
\centering
    \begin{tabular}{c}
    \includegraphics[clip, width=0.5 \textwidth]{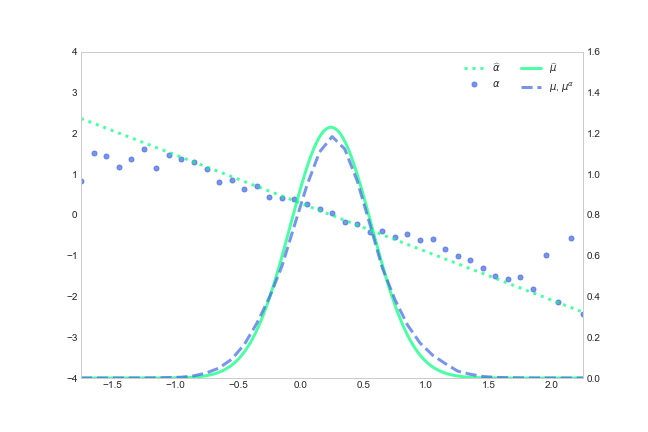}\\
    (a) $\rho_{k}^{\mu}<\rho_{k}^{Q}<\rho_{k}^{\mu^{\alpha}}$\\
    \begin{tabular}{cc}
    \includegraphics[clip, width=0.5 \textwidth]{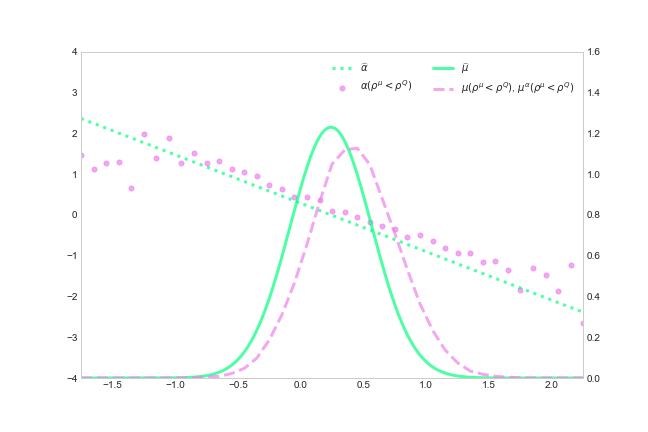}&
    \includegraphics[clip, width=0.5 \textwidth]{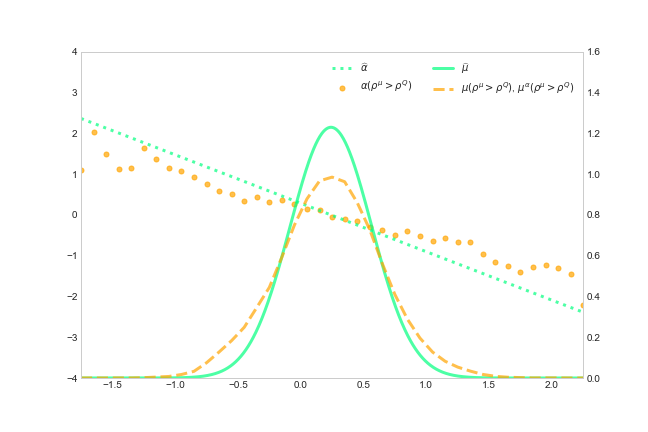}\\
    (b) $\rho_{k}^{\mu}=\rho_{k}^{\mu^{\alpha}}<\rho_{k}^{Q}$ & (c) $\rho_{k}^{\mu}=\rho_{k}^{\mu^{\alpha}}>\rho_{k}^{Q}$
    \end{tabular}
    \end{tabular}
    \caption{
    Control and distributions for the benchmark asymptotic MFCG learned by Algorithm~\ref{alg:rl_infinite}. The x-axis shows the state variable $x$, the left y-axis refers to the value of the control $\alpha(x)$, and the right y-axis marks the probability mass of $\mu(x)$ and  $\mu^{\alpha}(x)$. The green dotted line (labeled by $\hat{\alpha}$) is the theoretical control function and green curve (labeled by $\hat{\mu}$) shows the theoretical distribution of state, where the global distribution equals to the local distribution. The dots (labeled by $\alpha$) are the learned controls and the overlapping dashed curves (labeled by $\mu$ and $\mu^{\alpha}$) refer to the overlapping empirical global and local distributions learned by the algorithm, colored in blue, violet, and orange according to the selection of learning rates.
    }
    \label{fig:asym}
\end{figure}

\textbf{Figure~\ref{fig:asym}} is generated with $5$ runs of $K=\num{100000}$ episodes with the above setting. We report the average of the learned control and distribution in the last $\num{10000}$ episodes over the $5$ runs. In \textbf{Figure~\ref{fig:asym}}-(a) the control and distribution learned with the correct three-timescale rates by the algorithm are plotted against the theoretical optimal control obtained in Equation~\eqref{eq:alpha_sol_asym} and with the theoretical distribution. The control learned by the algorithm (blue dots) lies well along the theoretical control function (dotted green line) except for states never visited by the algorithm. That is, within the support of the distribution the algorithm learns the optimal control well. Also, we observe that the learned global distribution overlaps the local distribution (both dashed blue curve), and that both match the theoretical distribution (solid green curve). Therefore, the algorithm successfully learns the correct local and global distributions as well. \textbf{Figure~\ref{fig:asym}}-(b) and -(c) illustrate the failure of the algorithm in cases where the learning rates are misspecified as described in the last paragraph.

\subsection{Trader's Problem}\label{sec:traderPB}

As we did for the infinite horizon case, we consider a finite horizon extended game in continuous time and spaces. In particular, 
we reassess the renowned trader's execution problem presented in \cite{carmona2015probabilistic} starting with the finite number of players case. Instead of a game among traders, we consider a game between groups of traders as follows. Suppose there are $M$ homogeneous trading groups, each with $N$ traders trading on a single stock. Let the index tuple $(m,n)$ denotes the $n$-th trader in the $m$-th group.

Trader $(m,n)$ is controlling the drift term in the dynamic of her stock inventory,
\[
\mathrm{d}X_{t}^{m,n}=\alpha_{t}^{m}(X_{t}^{m,n})\mathrm{d}t+\sigma_{t}^{m,n}\mathrm{d}W_{t}^{m,n},
\]
by the trading rate $\alpha_{t}^{m}(X_{t}^{m,n})$ with which all the traders in the $m$-th group comply. Her cash position
$K_{t}^{m,n}$ evolves as
\[
\mathrm{d}K_{t}^{m,n}=-\left[\alpha_{t}^{m}(X_{t}^{m,n})S_{t}+c_{\alpha}\left(\alpha_{t}^{m}(X_{t}^{m,n})\right)\right]\mathrm{d}t,
\]
with $\alpha \mapsto c_{\alpha}(\alpha)$ a non-negative convex function representing the cost of trading at rate $\alpha$. The stock price is impacted by a function $h(\cdot)$ of the transactions
and follows the dynamic
\[
\mathrm{d}S_{t}=\frac{1}{M}\frac{1}{N}\sum_{i=1}^{M}\sum_{j=1}^{N}h\left(\alpha_{t}^{i}(X_{t}^{i,j})\right)\mathrm{d}t+\sigma_{t}^{0}\mathrm{d}W_{t}^{0}.
\]
The total wealth $V_{t}^{m,n}$ of her self-financing portfolio consists of her cash position and her stock value,
\[
V_{t}^{m,n}=K_{t}^{m,n}+X_{t}^{m,n}S_{t},
\]
with dynamic
\begin{align*}
\mathrm{d}V_{t}^{m,n} & =\mathrm{d}K_{t}^{m,n}+S_{t}\mathrm{d}X_{t}^{m,n}+X_{t}^{m,n}\mathrm{d}S_{t}\\
 & =\text{\ensuremath{\left[-c_{\alpha}\left(\alpha_{t}^{m}(X_{t}^{m,n})\right)+X_{t}^{m,n}\frac{1}{M}\frac{1}{N}\sum_{i=1}^{M}\sum_{j=1}^{N}h\left(\alpha_{t}^{i}(X_{t}^{i,j})\right)\right]}}\mathrm{d}t+S_{t}\sigma_{t}^{m,n}\mathrm{d}W_{t}^{m,n}+X_{t}^{m,n}\sigma_t^0\mathrm{d}W_{t}^{0}.
\end{align*}
We assume that the individual trader is subject to a running liquidation constraint modeled by a function $c_X$ of the average
shares held by her own group.
In this model, the individual trader's objective function is given by
\begin{align*}
&J^{m,n}\left(\alpha^{1},\dots,\alpha^{M}\right)\\
& =\mathbb{E}\left\{ \int_{0}^{T}c_{X}\left(\frac{1}{N}\sum_{j=1}^N X_{t}^{m,j}\right)\mathrm{d}t+g\left(X_{T}^{m,n}\right)-V_{T}^{m,n}\right\} \\
 & =\mathbb{E}\left\{ \int_{0}^{T}\left[c_{X}\left(\frac{1}{N}\sum_{j=1}^N X_{t}^{m,j}\right)+c_{\alpha}\left(\alpha_{t}^{m}(X_{t}^{m,n})\right)-X_{t}^{m,n}\frac{1}{M}\frac{1}{N}\sum_{i=1}^{M}\sum_{j=1}^{N}h\left(\alpha_{t}^{i}(X_{t}^{i,j})\right)\right]\mathrm{d}t+g\left(X_{T}^{m,n}\right)\right\} .
\end{align*}
In the limit of a large number of large groups without precising the relation between $M$ and $N$ (see Appendix \ref{sec:approx} for more details about this type of limit for a finite horizon linear-quadratic model), and assuming $\sigma_t^{m,n}=\sigma$, this problem leads to the following mixture of MFCG  problems: Minimize
\begin{align*}
J\left(\alpha; \theta\right)
 & =\mathbb{E}\left\{ \int_{0}^{T}\left[c_{X}\left(m^{\alpha, \theta}_t\right)+c_{\alpha}\left(\alpha_{t}\right)-X^{\alpha, \theta}_{t}\int h(a)d\theta_t(a)\right]\mathrm{d}t+g\left(X^{\alpha, \theta}_{T}\right)\right\} ,
\end{align*}
where $\theta_t$ is the law of the control $\alpha_t$,  $m_t^{\alpha, \theta} =\mathbb{E}(X_t^{\alpha, \theta})$, and
\begin{align*}
dX^{\alpha, \theta}_t&=\alpha_t dt+\sigma dW_t, \quad t\leq T,\quad X^{\alpha, \theta}_0=x.
\end{align*}
Note that the problem is of "MFG style" in control through $\theta_t$ and "MFC style" in state through $m_t^{\alpha, \theta}$.

In what follows we focus on the Linear-Quadratic case where $c_x(m)=\frac{c_X}{2}m^2$, $c_\alpha(\alpha)=\frac{c_\alpha}{2}\alpha^2$, $h(a)=c_h a$, and $g(x)=\frac{c_g}{2}x^2$, so that:
\begin{align*}
J\left(\alpha; \theta \right)
 & =\mathbb{E}\left\{ \int_{0}^{T}\left[\frac {c_{X}}{2}(m^{\alpha, \theta}_t)^2+\frac{c_{\alpha}}{2}\alpha_{t}^2-c_h X^{\alpha, \theta}_{t}\int ad\theta_t(a)\right]\mathrm{d}t+\frac{c_g}{2}(X^{\alpha, \theta}_{T})^2\right\} .
\end{align*}

\subsubsection{Results}

We consider the trader's problem with the choice of parameters: $c_\alpha=1$, $c_X=0.75$, $c_h=1.25$, $c_g=1$, and with a volatility for the state dynamic $\sigma = 0.75$. We test three distributions for the initial inventory $X_0$: Gaussian with mean $x_0 = 0$, $0.5$, and $1$ and the same standard deviation $\sigma = 0.5$. The terminal time is $T=1$, and we choose a time grid $\tau = \{0, \Delta t, \dots, T\}$ with time step $\Delta t = 1/16$. We discretize the state space into $\X = \{x_0=-2, \dots, x_{|{\X}|-1}=2.5\}$, and the action space into $\A = \{a_0 = -2, \dots, a_{|{\A}|-1}=1.5\}$, where the step sizes are $\Delta x = \Delta a = \sqrt{\Delta t}=1/4$. The triplet of the learning rates is chosen as $(\omega^{\theta}, \omega^{Q}, \omega^{\mu}) = (0.85, 0.55, 0.15)$. For the $\epsilon$-greedy policy we choose $\epsilon=0.05$. We run the experiment $10$ times, each with $K=\num{200000}$ episodes. We average the control and state distributions learned by Algorithm~\ref{alg:rl_finite} over the last $\num{10000}$ episodes and over $10$ runs. We report the results in \textbf{Figure~\ref{fig:trader}}. We present the results for every time step in $\tau$, except for the last time step $T$. The subplots are ordered by time, from left to right and top to bottom. Note that the theoretical optimal control $\hat{\alpha}_t$  (dotted green line) changes over time. As time increases, the slope and intercept of $\hat{\alpha}_t$ increase. Also, the theoretical local state distribution $\mu_t$ (green curve) under the optimal control changes over time. As time increases from $0$ to $T$, the center of $\mu_t$ moves towards zero and the standard deviation increases. To evaluate the effectiveness of Algorithm~\ref{alg:rl_finite}, we compare the learned action (blue dots) and the learned local state distribution (dashed blue curve) with their theoretical counterparts. Again we observe that except for the tails of the distribution, the control learned by the algorithm is very close to the theoretical value. This means that the algorithm successfully learns the optimal control for states that are frequently sampled. Also, we see that the dashed blue curve perfectly overlaps with the solid green curve, hence the algorithm succeeds in capturing the evolution of the state distribution under the correctly learned control.

\begin{figure}[H]
    \centering
    \begin{subfigure}{0.7\textwidth}
        \includegraphics[trim= 100 100 100 100, width=\linewidth]{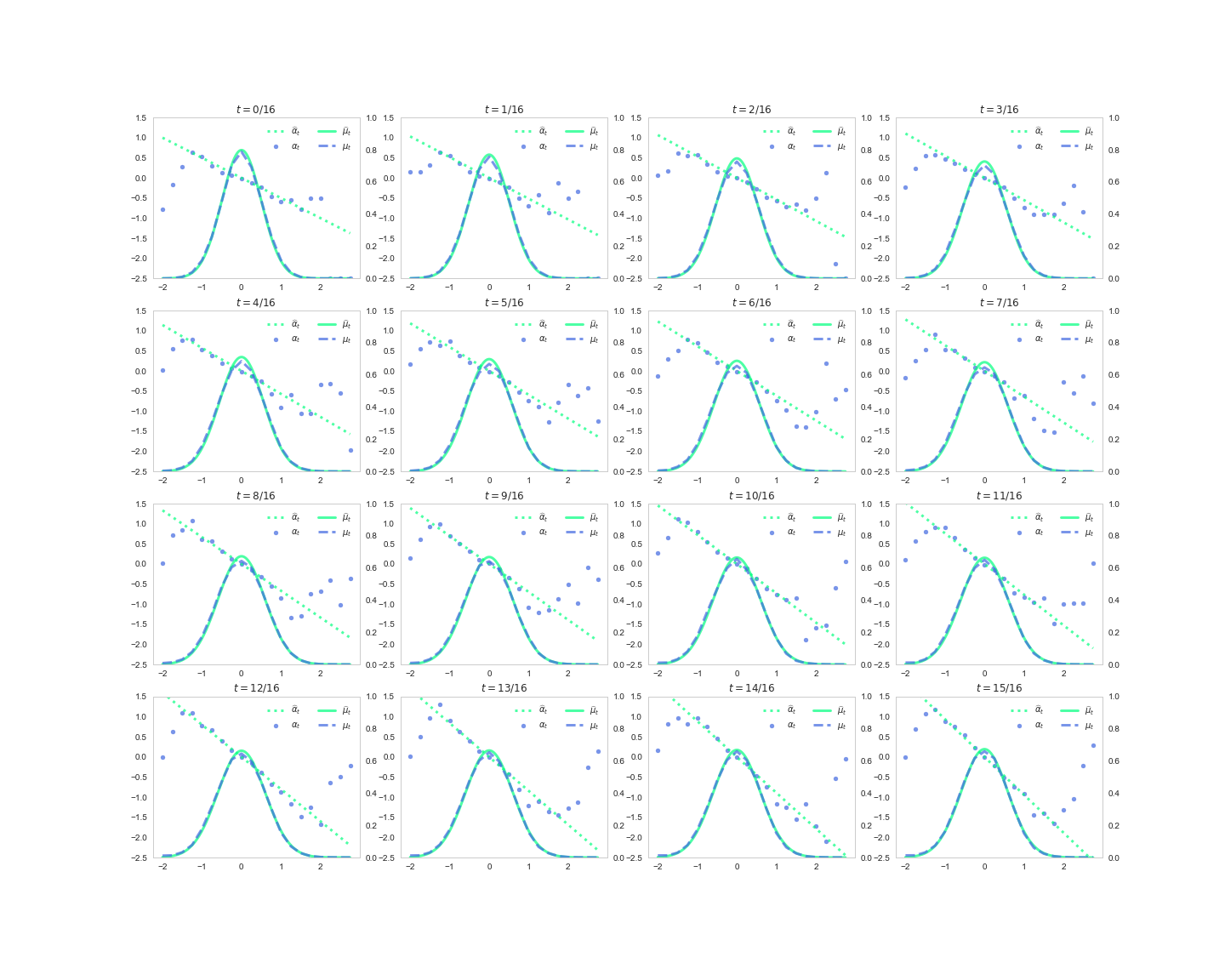}
        \subcaption{$x_0$ = 0}
    \end{subfigure}
\end{figure}
\begin{figure}[H] \ContinuedFloat
    \centering
    \begin{subfigure}{0.7\textwidth}
        \includegraphics[trim= 100 100 100 100, width=\linewidth]{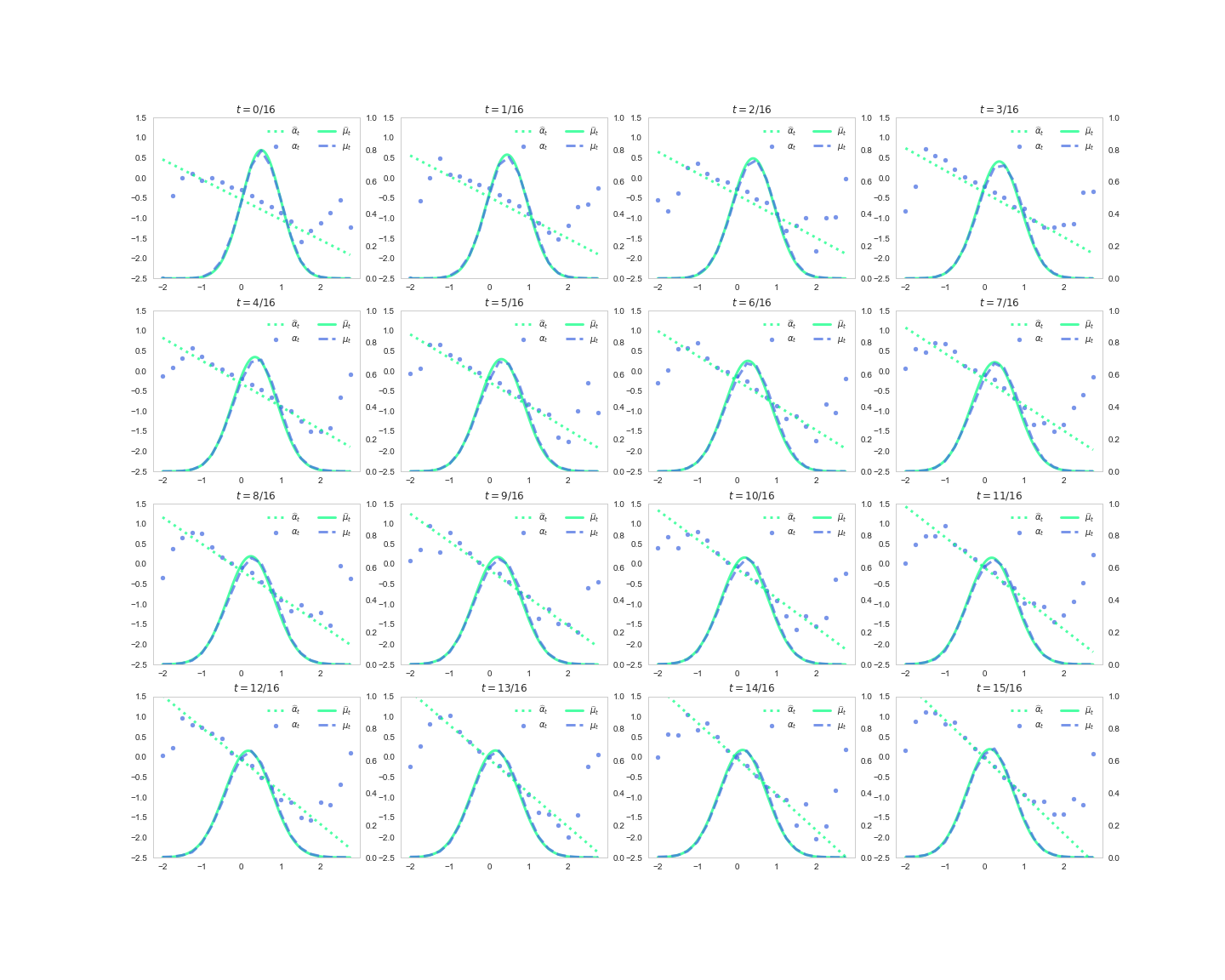}
        \subcaption{$x_0$ = 0.5}
    \end{subfigure}
\end{figure}
\begin{figure}[H] \ContinuedFloat
    \centering
    \begin{subfigure}{0.7\textwidth}
        \includegraphics[trim= 100 100 100 100, width=\linewidth]{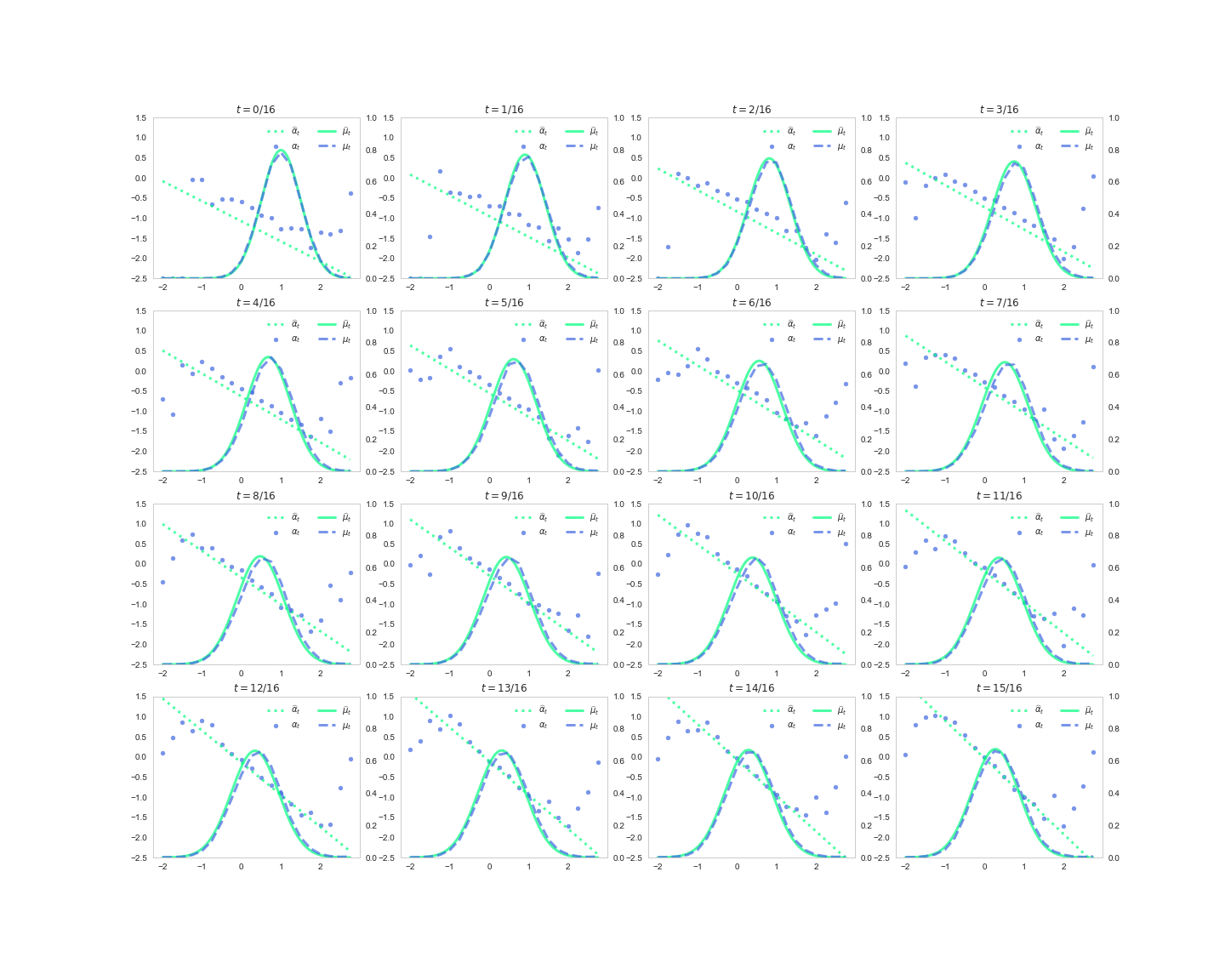}
        \subcaption{$x_0$ = 1}
    \end{subfigure}
    \caption{Control and distributions for the trader's MFCG learned by Algorithm~\ref{alg:rl_finite}. The x-axis shows the state variable $x$, the left y-axis refers to the value of the control $\alpha(x)$, and the right y-axis marks the probability mass of state, $\mu(x)$. The dotted green lines (labeled by $\hat{\alpha}_t$) are the theoretical control function and the blue dots (labeled by $\alpha_t$) are the learned control. The green curves (labeled by $\hat{\mu}_t$) show the theoretical distributions of state and dashed blue curves (labeled by $\mu_t$) refer to the empirical distribution of state learned by the algorithm.}
    \label{fig:trader}
\end{figure}

\section{Conclusion}

We have introduced a type of mean field control games (MFCG) that models competitive games between a large number of large collaborative groups. It turns out that the two-timescale reinforcement learning algorithm (U2-MF-QL) that was proposed in \cite{AFL2021} for infinite horizon problems and in \cite{AFL2022} for finite horizon extended problems, for learning either MFG or MFC problems, is naturally adapted for learning MFCG problems by managing three learning rates in the three-timescale reinforcement learning algorithm (U3-MF-QL) proposed in this paper. We illustrate the results with linear quadratic problems for which we have explicit formulas. In particular, a new type of trader problem is presented. The theory associated for MFCGs is a work in progress \cite{MultiTimescaleConvergenceproof}, as well as an actor-critic version of the U3-MF-QL algorithm in the context of continuous spaces \cite{ADFLcontinuous}.

\bibliographystyle{chicago}
\bibliography{bibtex}

\appendix

\section{Linear-Quadratic Example}\label{appendix}
In this Appendix we provide additional details about the LQ example presented in Section \ref{sec:motiv}, here in a finite horizon setting. In particular, we explain the relation between the finite-player game and its corresponding MFCG limiting problem.

\subsection{The Finite-Player Model}\label{sec:finite}
There are $M$ competitive groups each made of $N$ collaborative players ($m=1,\cdots,M$ is the group number and $n=1,\cdots,N$ is the player number within the group). The dynamics of the state of player $(m,n)$ is:
\[
dX^{m,n}_t=\alpha^{m,n}_tdt +\sigma dW^{m,n}_t, \quad X^{m,n}_0 \sim \mu_0,
\]
where we aim at an open-loop equilibrium (the $\alpha$'s are adapted to the $W$'s). The objective of the collaborative group $m$ is to minimize
\[
J^m(\mathbf\alpha)=\frac{1}{N}\sum_{n=1}^N\E\int_0^T\left\{\frac{1}{2}(\alpha^{m,n}_t)^2+\frac{c_1}{2}(\bar\mu_t-X^{m,n}_t)^2+\frac{c_2}{2}(\bar\mu^m_t)^2\right\}dt
\]
where $\bar\mu^m_t=\frac{1}{N}\sum_{n=1}^NX^{m,n}_t$ is the empirical mean of group $m$, and $\bar\mu_t=\frac{1}{MN}\sum_{m=1}^M\sum_{n=1}^NX^{m,n}_t$ is the empirical mean of the total population ($c_1$ and $c_2$ are positive constants and we assume a zero terminal condition for simplicity).

Accordingly, we introduce the Hamiltonian $H^m$ of group $m$:
\[
H^{m,n}=\sum_{m'=1}^M\sum_{n'=1}^N\alpha^{m',n'}y^{m,m',n,n'}+
\sum_{n=1}^N\frac{1}{2}(\alpha^{m,n})^2+\frac{c_1}{2}\sum_{n=1}^N(\bar\mu-x^{m,n})^2+\frac{c_2}{2}N(\bar\mu^m)^2,
\]
where $y^{m,m',n,n'}$ are the adjoint variables.

Minimizing the Hamiltonian $H^{m,n}$ with respect to $\alpha^{m,n}$, we get
\[
\frac{\partial H^{m,n}}{\partial \alpha^{m,n}}=y^{m,m,n,n}+\alpha^{m,n}=0 \implies \hat\alpha^{m,n}=-y^{m,m,n,n}.
\]
The Backward Stochastic Differential Equation (BSDE) that the  adjoint process $Y^{m,m',n,n'}_t$ must satisfy is:
\begin{align*}
dY^{m,m',n,n'}_t&=-\partial_{x^{m',n'}}H^{m,n} dt+\sum_{m"=1}^M\sum_{n"=1}^NZ^{m,m',n,n',m",n"}_tdW^{m",n"}_t\\
&=-\left[c_1\sum_{k=1}^N(\bar\mu_t-X_t^{m,k})(\frac{1}{MN}-\delta_{\{m'=m, k=n'=n\}})+c_2N\bar\mu^m_t
(\frac{1}{N}\delta_{\{m'=m,n'=n\}})
\right]dt\\
&+\sum_{m"=1}^M\sum_{n"=1}^NZ^{m,m',n,n',m",n"}_tdW^{m",n"}_t,
\end{align*}
with a zero terminal condition $Y^{m,m',n,n'}_T=0$. The $Z$-processes are part of the solution and must be adapted.

The diagonal adjoint process $Y^{m,m,n,n}_t\equiv Y^{m,n}_t$  satisfies
\begin{align}
dY^{m,n}_t&=-\left[c_1\sum_{k=1}^N(\bar\mu_t-X_t^{m,k})(\frac{1}{MN}-\delta_{\{k=n\}})+c_2\bar\mu^m_t
\right]dt +\sum_{m"=1}^M\sum_{n"=1}^NZ^{m,m,n,n,m",n"}_tdW^{m",n"}_t \nonumber\\
&=\left[c_1\left((1-\frac{1}{M})\bar\mu_t+\frac{1}{M}\bar\mu^m_t-X^{m,n}_t\right)-c_2\bar\mu^m_t\right]dt +\sum_{m"=1}^M\sum_{n"=1}^NZ^{m,m,n,n,m",n"}_tdW^{m",n"}_t\label{eqY}
\end{align}
We omit the non-diagonal adjoint processes which can be treated analogously, and we formulate the 
ansatz:
\[
Y^{m,n}_t=\eta_tX^{m,n}_t+\phi_t\bar\mu_t+\xi_t\bar\mu^m_t,
\]
where $\eta_t$, $\phi_t$ and $\xi_t$ are deterministic functions to be determined.
We have
\begin{align*}
d\bar\mu^m_t&=\frac{1}{N}\sum_{n=1}^NdX^{m,n}_t
=-\bar Y^m_tdt+\frac{\sigma}{N}\sum_{n=1}^NdW^{m,n}_t
\end{align*}
where
\[
\bar Y^m_t=\frac{1}{N}\sum_{n=1}^NY^{m,n}_t=\phi_t\bar\mu_t+(\eta_t+ \xi_t)\bar\mu^m_t,
\]
and
\begin{align*}
d\bar\mu_t&=\frac{1}{MN}\sum_{n=1}^N\sum_{m=1}^MdX^{m,n}_t
=-\bar{Y}_t dt+\frac{\sigma}{MN}\sum_{n=1}^N\sum_{m=1}^MdW^{m,n}_t,
\end{align*}
where
\[
\bar{Y}_t=\frac{1}{MN}\sum_{n=1}^N\sum_{m=1}^MY^{m,n}_t=\left(\eta_t+\phi_t+\xi_t\right)\bar\mu_t
\]
Differentiating the ansatz gives:
\begin{align*}
dY^{m,n}_t&=\eta'_tX^{m,n}_tdt+\eta_tdX^{m.n}_t+\phi'_t\bar\mu_tdt+\phi_td \bar\mu_t+ \xi'_t\bar\mu^m_tdt+\xi_td\bar\mu^m_t\\
&=\eta'_tX^{m,n}_tdt-\eta_t(\eta_tX^{m,n}_t+\phi_t\bar\mu_t+\xi_t\bar\mu^m_t)dt+\phi'_t\bar\mu_tdt -\phi_t\left(\eta_t+\phi_t
+\xi_t\right)\bar\mu_tdt\\
&+ \xi'_t\bar\mu^m_tdt-\xi_t(\phi_t\bar\mu_t+(\eta_t+ \xi_t)\bar\mu^m_t)dt
+d \mbox{Mart}\\
&=\left[\eta'_t-\eta_t^2\right]X^{m,n}_tdt+\left[\phi'_t-   \eta_t\phi_t -\phi_t\left(\eta_t+\phi_t
+\xi_t\right)-\xi_t\phi_t \right]\bar\mu_tdt+\left[\xi'_t-\eta_t\xi_t-\xi_t(\eta_t+\xi_t)\right]\bar\mu^m_tdt
+d \mbox{Mart}\\
&=\left[\eta'_t-\eta_t^2\right]X^{m,n}_tdt+\left[\phi'_t-   \phi_t^2 -2\eta_t\phi_t     -2\xi_t\phi_t \right]\bar\mu_tdt+
\left[\xi'_t-\xi_t^2 -2\eta_t\xi_t\right]\bar\mu^m_tdt+d \mbox{Mart}.
\end{align*}
Comparing the drift terms with the previous expression \eqref{eqY} for $dY^{m,n}_t$, we get the following system of Riccati equations for which explicit solutions can be obtained (omitted here):
\begin{align*}
\eta'_t-\eta_t^2&=-c_1, \quad \eta_T=0,\\
\phi'_t-   \phi_t^2 -2\phi_t(\eta_t+\phi_t )   &=-c_2+\frac{1}{M}c_1,\quad \phi_T=0,\\
\xi'_t-\xi_t^2 -2\eta_t\xi_t&= c_1(1-\frac{1}{M}),\quad \xi_T=0.
\end{align*}
As usual the $Z$'s processes are deterministic (hence adapted) and identified by matching the martingale terms. 

One can define $\zeta_t=\eta_t+\xi_t$, so that the system of ODEs becomes
\begin{align}
\eta'_t-\eta_t^2&=-c_1, \quad \eta_T=0,\nonumber\\
\phi'_t-   \phi_t^2 -2\phi_t\zeta_t   &=-c_2+\frac{1}{M}c_1,\quad \phi_T=0,\label{ODE1}\\
\zeta'_t-\zeta_t^2 &= -\frac{1}{M}c_1,\quad \xi_T=0,\nonumber
\end{align}
which highlights the limit $M\to\infty$ where $\zeta_t$ vanishes.

\subsection{The Corresponding Limiting MFCG}\label{sec:limit}
To the previous finite-player model, we propose to associate the following MFCG problem:
For a fixed flow of distributions $\mu= (\mu_t)$, one agent controls her state given by
\[
dX^{\alpha, \mu}_t=\alpha_tdt+\sigma dW_t, \quad X^{\alpha, \mu}_0\sim \mu_0.
\]
The agent solves the MKV control problem which consists in minimizing
\[
J(\alpha; \mu)=\E\int_0^T\left\{ \frac{1}{2}\alpha^2_t+\frac{c_1}{2}(\bar\mu_t-X^{\alpha, \mu}_t)^2+\frac{c_2}{2}(\E(X^{\alpha, \mu}_t))^2\right\}dt
\]
where $\bar\mu_t=\int x\mu(dx)$. One then solves the fixed point condition:
\[
\E(X^{\hat\alpha, \mu}_t)=\bar\mu_t ,\quad \forall t\leq T.
\]
Introducing the adjoint process $Y_t$ and using the lighter notation $X_t$ for the state process, the optimal strategy $\hat\alpha_t$ is given by $-Y_t$ which satisfies the BSDE
\[
dY_t=\left[c_1(\E(X_t)-X_t)-c_2\E(X_t)\right]dt +Z_tdW_t, \quad Y_T=0
\]
The term $c_2\E(X_t)$ comes form the differentiation of $\frac{c_2}{2}(\E(X^{\alpha, \mu}_t))^2$ with respect to the measure.

One verifies easily that the solution is 
\[
Y_t=-\eta_t
\left(\bar\mu_t-X_t\right)+\phi_t\bar\mu_t
\]
with 
\begin{align}
\eta'_t-\eta_t^2&=-c_1, \quad \eta_T=0,\\
\phi'_t-   \phi_t^2  &=-c_2,\quad \phi_T=0,\label{ODE2}
\end{align}
and 
\[
d\bar\mu_t=-\phi_t \bar\mu_t,\quad \bar\mu_0=x_0,
\]
that is
\[
\bar\mu_t=x_0e^{-\int_0^t\phi_sds}.
\]
Note that the functions $\eta$ and $\phi$ are given explicitly by
\[
\eta_t=\sqrt{c_1}\frac{e^{2\sqrt{c_1}(T-t)}-1}{e^{2\sqrt{c_1}(T-t)}+1},\quad \phi_t=\sqrt{c_2}\frac{e^{2\sqrt{c_2}(T-t)}-1}{e^{2\sqrt{c_2}(T-t)}+1}.
\]

\subsection{From Finite-Player to MFCG}\label{sec:approx}
The limit $N\to\infty$ ensures that $\bar\mu^m_t=\bar\mu_t$ for every $m$, and the limit $M\to \infty$ ensures that the coefficient functions  given by \eqref{ODE1} converge to those given by \eqref{ODE2}.

Our goal is to show that the strategy obtained from the limiting problem in Section \ref{sec:limit} provides an $\epsilon$-Nash equilibrium for the finite-player game described in Section \ref{sec:finite}.

Here we use the notation $(\eta^\infty_t,\phi^\infty_t, \bar\mu^\infty_t)$ for the quantities obtained in the system of equations \eqref{ODE2} (not to be confused with the corresponding quantities obtained in \eqref{ODE1}.

We denote by $\alpha^\infty$ the optimal strategy obtained in Section \ref{sec:limit}, that is: 
\[
\alpha^\infty_t=-Y_t=\eta^\infty_t
\left(\bar\mu^\infty_t-X_t\right)-\phi^\infty_t\bar\mu^\infty_t,
\]
which we apply to all the players in the finite-player game. The value function for the $m$-th group is given by
\[
J^{m}(\alpha^\infty)=\frac{1}{N}\sum_{n=1}^N\E\int_0^T\left\{\frac{1}{2}(\alpha^{m,n}_t)^2+\frac{c_1}{2}(\bar\mu_t-X^{m,n}_t)^2+\frac{c_2}{2}(\bar\mu^m_t)^2\right\}dt,
\]
where 
\[
\alpha^{m,n}_t=\eta^\infty_t
\left(\bar\mu^\infty_t-X^{m,n}_t\right)-\phi^\infty_t\bar\mu^\infty_t,
\]
 and
\[
dX^{m,n}_t=\alpha^{m,n}_tdt +\sigma dW^{m,n}_t, \quad X^{m,n}_0 \sim \mu_0.
\]
Note that $\bar\mu^m_t=\frac{1}{N}\sum_{n=1}^NX^{m,n}_t$ is given by
\[
d\bar\mu^m_t=[\eta^\infty_t
\left(\bar\mu_t-\bar\mu^m_t\right)-\phi^\infty_t\bar\mu_t]dt+\frac{1}{N}\sum_{n=1}^NdW^{m,n}_t, \quad \bar\mu_0=\frac{1}{N}\sum_{n=1}^NX^{m,n}_0,
\]
and $\bar\mu_t=\frac{1}{MN}\sum_{m=1}^M\sum_{n=1}^NX^{m,n}_t$ is given by
\[
d\bar\mu_t=[-\phi^\infty_t\bar\mu_t]dt+\frac{1}{MN}\sum_{m=1}^M\sum_{n=1}^NdW^{m,n}_t, \quad \bar\mu_0=\frac{1}{MN}\sum_{m=1}^M\sum_{n=1}^NX^{m,n}_0.
\]
Now we consider a strategy $(\alpha^{\infty,-m},\beta^m)$ where the players from group $m$ use $\beta^{m,n}_t$ instead of $\alpha^{m,n}_t$, and the players from the other groups $m'$ continue using $\alpha^{m',n}_t$. 
We denote by $\tilde{X}^{m,n}_t$ the state of player $(m,n)$ which satisfies
\[
d\tilde{X}^{m,n}_t=\beta^{m,n}_tdt +\sigma dW^{m,n}_t, \quad X^{m,n}_0 \sim \mu_0.
\]
We denote the corresponding group empirical mean by 
$\tilde\mu^m_t$, and the population empirical mean by $\tilde\mu_t$. We also denote:
\[
f(x,\alpha,\bar\mu,\bar\mu^m)=\frac{1}{2}\alpha^2+ F(x,\bar\mu,\bar\mu^m),\quad F(x,\bar\mu,\bar\mu^m)=\frac{1}{2}
\left\{c_1(\bar\mu-x)^2+c_2(\bar\mu^m)^2\right\}
\]
\underline{Principle of the proof:}

First, show that for $\epsilon>0$ there exists $M_0$ and $N_0$ such that for $M\geq M_0$ and $N\geq N_0$, we have 
\begin{equation}\label{ineq1}
\left |
\frac{1}{N}\sum_{n=1}^N\E\int_0^T\left(F(\tilde{X}^{m,n}_t,\tilde\mu_t,\tilde\mu^m_t)-
F(X^{m,n}_t,\bar\mu_t,\bar\mu_t^m)\right)dt
\right |<\frac{\epsilon}{2}
\end{equation}
where in the first term the strategy $(\alpha^{\infty,-m},\beta^m)$ is used, while in the second term the strategy $(\hat\alpha^{-m},\beta^m)$ is used with $\hat{\alpha}$ the optimal strategy obtained in Section \ref{sec:finite} for the finite-player game. Adding $\frac{1}{2}(\beta^{m,n}_t)^2$ to both terms, we obtain:
\[
J^{m}(\alpha^{\infty,-m},\beta^m)> J^{m}(\hat\alpha^{-m},\beta^m)-\frac{\epsilon}{2}
\]
Using the fact that $\hat\alpha$ is a Nash equilibrium for the finite-player game, we get
\[
J^{m}(\alpha^{\infty,-m},\beta^m)> J^{m}(\hat\alpha)-\frac{\epsilon}{2}
\]
As in the first step we can derive
\begin{equation}\label{ineq2}
\left | J^{m}(\alpha^{\infty})- J^{m}(\hat\alpha)\right | <\frac{\epsilon}{2}
\end{equation}
and, therefore
\[
J^{m}(\alpha^{\infty,-m},\beta^m)> J^{m}(\alpha^\infty)-\epsilon
\]
that is $\alpha^\infty$ is an $\epsilon$-Nash equilibrium for the finite-player game. 

Of course, \eqref{ineq1} and \eqref{ineq2} require some technical work which will be given in a general setting in \cite{MultiTimescaleConvergenceproof}.

Finally we observe that the limits $N\to\infty$ and $M\to \infty$ can be taken sequentially.

If $N\to\infty$ for $M$ fixed, we obtain a game between competitive MKV agents called \textit{Mean Field Type Game} in \cite{tembine}.
Then, our limit describes the MFG limit between these MKV agents.

If $N$ is fixed and $M\to\infty$, on can consider each group as one player in  the higher dimension $N$ and this is a classical MFG. Our MFCG describes the subsequent limit $N\to\infty$.

\section{Analytic solutions}\label{sec:appendix_solutions}
In this Appendix we provide details about the solutions of the problems discussed in Sections \ref{sec: Numerical results}.

\subsection{Solution of the Asymptotic Problem}\label{appendix:LQ}

The corresponding HJB equation is given by 
\[
\beta V(x)-H(x,\alpha,\mu,\mu^{\alpha, \mu})-\int_{\mathds{R}}\frac{\partial H}{\partial\mu^{\alpha, \mu}}H\left(h,\alpha, \mu,\mu^{\alpha, \mu}\right)(x)\mathrm{d}\mu^{\alpha, \mu}(h)=0,
\]
with the Hamiltonian
\begin{align*}
H(x,\alpha, \mu,\mu^{\alpha, \mu}) & =\inf_{\alpha}\left\{ \mathcal{A}^{X}V(x)+f(x,\alpha,\mu,\mu^{\alpha, \mu})\right\} \\
 & =\inf_{\alpha}\left\{ \alpha\dot{V}(x)+\frac{1}{2}\sigma^{2}\ddot{V}(x)+\frac{1}{2}   \alpha^{2}+c_{1}(x-c_{2}m)^{2}+c_{3}(x-c_{4})^{2}+\tilde{c}_{1}(x-\tilde{c}_{2}{m}^{\alpha, \mu})^{2}+\tilde{c}_{5}\left({m}^{{\alpha, \mu}}\right)^{2}\right\} \\
 & =-\frac{1}{2}\dot{V}(x)^{2}+\frac{1}{2}\sigma^{2}\ddot{V}(x)+c_{1}(x-c_{2}m)^{2}+c_{3}(x-c_{4})^{2}+\tilde{c}_{1}(x-\tilde{c}_{2}{m}^{\alpha, \mu})^{2}+\tilde{c}_{5}\left({m}^{{\alpha, \mu}}\right)^{2},
\end{align*}
and the derivative with respect to $\mu^{\alpha, \mu}$ due to the MFC part, calculated at the optimal $\hat\alpha(x)=-\dot{V}(x)$ as follows:
\begin{align*}
\frac{\partial H}{\partial\mu^{\alpha, \mu}}\left(h,-\dot{V}(h),\mu,\mu^{\alpha, \mu}\right) & =\frac{\partial}{\partial\mu^{\alpha, \mu}}\left(\tilde{c}_{1}(h-\tilde{c}_{2}{m}^{\alpha, \mu})^{2}+\tilde{c}_{5}\left({m}^{\alpha, \mu}\right)^{2}\right)(x)\\
 & =\frac{\partial}{\partial\mu^{\alpha, \mu}}\left(\tilde{c}_{1}\left(h-\tilde{c}_{2}\int_{\mathds{R}}y\mathrm{d}\mu^{\alpha, \mu}(y)\right)^{2}+\tilde{c}_{5}\left(\int_{\mathds{R}}y\mathrm{d}\mu^{\alpha, \mu}(y)\right)^{2}\right)(x)\\
 & =-2\tilde{c}_{1}\tilde{c}_{2}x\left(h-\tilde{c}_{2}\int_{\mathds{R}}y\mathrm{d}\mu^{\alpha, \mu}(y)\right)+2\tilde{c}_{5}x\int_{\mathds{R}}y\mathrm{d}\mu^{\alpha, \mu}(y)\\
 & =-2\tilde{c}_{1}\tilde{c}_{2}x\left(h-\tilde{c}_{2}{m}^{\alpha, \mu}\right)+2\tilde{c}_{5}x{m}^{\alpha, \mu},
\end{align*}
and
\begin{align*}
\int_{\mathds{R}}\frac{\partial H}{\partial\mu^{\alpha, \mu}}\left(h,-\dot{V}(h),\mu,\mu^{\alpha, \mu}\right)(x)\mathrm{d}\mu^{\alpha, \mu}(h) & =-2\tilde{c}_{1}\tilde{c}_{2}(1-\tilde{c}_{2})x{m}^{\alpha, \mu}+2\tilde{c}_{5}x{m}^{\alpha, \mu},
\end{align*}
Finally, the HJB equation reduces to:
\begin{equation}
\beta V(x)+\frac{1}{2}\dot{V}(x)^{2}-\frac{1}{2}\sigma^{2}\ddot{V}(x)-c_{1}(x-c_{2}m)^{2}-c_{3}(x-c_{4})^{2}-\tilde{c}_{1}(x-\tilde{c}_{2}m^{\alpha, \mu})^{2}-\tilde{c}_{5}\left({m}^{\alpha, \mu}\right)^{2}+2\tilde{c}_{1}\tilde{c}_{2}(1-\tilde{c}_{2})x{m}^{\alpha, \mu}-2\tilde{c}_{5}x{m}^{\alpha, \mu}=0.\label{eq:HJB}
\end{equation}
Using the following ansatz for the value function and its derivatives
\begin{align}
V(x) & =\Gamma_{2}x^{2}+\Gamma_{1}x+\Gamma_{0},\nonumber \\
\dot{V}(x) & =2\Gamma_{2}x+\Gamma_{1},\label{eq:value_ansatz}\\
\ddot{V}(x) & =2\Gamma_{2},\nonumber 
\end{align}
we obtain
 the optimal control 
\begin{equation}
\hat\alpha(x)=-\dot{V}(x)=-2\Gamma_{2}x-\Gamma_{1},\label{eq:alpha_optimal}.
\end{equation}
Plugging the ansatz \eqref{eq:value_ansatz} into the HJB \eqref{eq:HJB}
we have 
\begin{align*}
\left(\beta\Gamma_{2}+2\Gamma_{2}^2 -(c_{1}+c_{3}+\tilde{c}_{1})\right)x^{2}\\
+\left(\beta\Gamma_{1}+2\Gamma_{2}\Gamma_{1}+2c_{1}c_{2}m+2\tilde{c}_{1}\tilde{c}_{2}{m}^{\alpha, \mu}+2c_{3}c_{4}+2\tilde{c}_{1}\tilde{c}_{2}(1-\tilde{c}_{2}){m}^{\alpha, \mu}-2\tilde{c}_{5}{m}^{\alpha, \mu}\right)x\\
+\beta\Gamma_{0}+\frac{1}{2}\Gamma_{1}^{2}-\sigma^{2}\Gamma_{2}-c_{1}c_{2}^{2}m^{2}-\left(\tilde{c}_{1}\tilde{c}_{2}^{2}+\tilde{c}_{5}\right)\left({m}^{\alpha, \mu}\right)^{2}-c_{3}c_{4}^{2} & =0.
\end{align*}
The solution is given by 
\begin{align*}
\Gamma_{2} & =\frac{-\beta+\sqrt{\beta^{2}+8\left(c_{1}+c_{3}+\tilde{c}_{1}\right)}}{4},\\
\Gamma_{1} & =\frac{2\tilde{c}_{5}{m}^{\alpha, \mu}-2\tilde{c}_{1}\tilde{c}_{2}(2-\tilde{c}_{2}){m}^{\alpha, \mu}-2c_{1}c_{2}m-2c_{3}c_{4}}{\beta+2\Gamma_{2}},\\
\Gamma_{0} & =\frac{c_{1}c_{2}^{2}m^{2}+\left(\tilde{c}_{1}\tilde{c}_{2}^{2}+\tilde{c}_{5}\right)\left({m}^{\alpha, \mu}\right)^{2}+\sigma^{2}\Gamma_{2}-\frac{1}{2}\Gamma_{1}^{2}+c_{3}c_{4}^{2}}{\beta}.
\end{align*}
Taking the expectation of the dynamics of $X_{t}^{\alpha, \mu}$
with the control $\hat\alpha(x)$,  we obtain the following ODE for ${m}^{\hat\alpha, \mu}$:
\[
\dot{{m}}_{t}^{\hat\alpha, \mu}=-2\Gamma_{2}{m}_{t}^{\hat\alpha, \mu}-\Gamma_{1},
\]
which is solved by 
\begin{align*}
{m}^{\hat\alpha, \mu} & =\lim_{t\rightarrow\infty}{m}_{t}^{\hat\alpha, \mu}=\lim_{t\rightarrow\infty}\left(-\frac{\Gamma_{1}}{2\Gamma_{2}}+\left({m}_{0}^{\hat\alpha, \mu}+\frac{\Gamma_{1}}{\Gamma_{2}}\right)\mathrm{e}^{-2\Gamma_{2}t}\right)\\
 & =-\frac{\Gamma_{1}}{2\Gamma_{2}}=-\frac{2\tilde{c}_{5}{m}^{\hat\alpha, \mu}-2\tilde{c}_{1}\tilde{c}_{2}(2-\tilde{c}_{2}){m}^{\hat\alpha, \mu}-2c_{1}c_{2}m-2c_{3}c_{4}}{2\Gamma_{2}(\beta+2\Gamma_{2})}.
\end{align*}
From  the fixed point condition $m={m}^{\hat\alpha, \mu}$,
we deduce
\[
\hat{m}={m}^{\hat\alpha, \hat\mu}=\frac{c_{3}c_{4}}{c_{1}(1-c_{2})+\tilde{c}_{1}(1-\tilde{c}_{2})^2+c_{3}+\tilde{c}_{5}},
\]
and the explicit form of the optimal control \eqref{eq:alpha_optimal}
\begin{align} \label{eq:alpha_sol_asym}
\hat\alpha(x)=-2\Gamma_{2}\left(x-\hat{m}\right).
\end{align}
Note that $\mu^{\hat\alpha, \hat\mu}={\mathcal N}\left(\hat{m}, \frac{\sigma^2}{4\Gamma_2}\right)$ is the limiting distribution of the OU process $(X_t^{\hat\alpha, \hat\mu})$.

\subsection{Solution of the Traders' Problem}\label{sec:appendix_trader}

In order to solve this problem, one first "freezes" the flow $(\theta_t)$ as in the MFG problem, and then solves the control problem which is of MKV type due to the term $m_t =\mathbb{E}(X_t)$ of MFC style.
Differentiating the corresponding Hamiltonian with respect to $\alpha$, one gets
\begin{align*}
\hat\alpha_t&=-\frac{1}{c_\alpha}Y_t.
\end{align*}
On the other hand,
\[
dY_t=-\left(-c_h \mathbb{E}[\hat\alpha_t]+c_X\mathbb{E}[X_t]\right)dt+Z_tdW_t,
\]
which leads to the following FBSDE:
\[
\begin{cases}
dX_t=-\frac{1}{c_\alpha} Y_tdt+\sigma dW_t, \quad X_0\sim \mu_0,\\
dY_t=-\left(\frac{c_h}{c_\alpha}\mathbb{E}[Y_t]
+c_X\mathbb{E}[X_t]\right)dt+Z_t dW_t,\quad  Y_T=c_g X_T.
\end{cases}
\]
Note that this is a different system than the one studied in \cite[Section 6.2.]{AFL2022}.
Taking expectation in this system one obtains:
\[
\begin{cases}
d\mathbb{E}[X_t]=-\frac{1}{c_\alpha}\mathbb{E}[Y_t]dt, \quad \mathbb{E}[X_0]=x_0,\\
d\mathbb{E}[Y_t]=-\left(\frac{c_h}{c_\alpha}\mathbb{E}[Y_t]
+c_X\mathbb{E}[X_t]\right)dt,\quad \mathbb{E}[Y_T]=c_g\mathbb{E}[ X_T].
\end{cases}
\]
Solving this system leads to
\[
\mathbb{E}[Y_t]=\bar\eta(t)\mathbb{E}[X_t],
\]
where
\begin{equation*}
\bar{\eta}_t= \frac{-C (e^{(\delta^+-\delta^-)(T-t)}-1)-c_g(\delta^+e^{(\delta^+-\delta^-)(T-t)}-\delta^-)}{(\delta^-e^{(\delta^+-\delta^-)(T-t)}-\delta^+)-c_gB(e^{(\delta^+-\delta^-)(T-t)}-1)},
\end{equation*}
for $t\in[0,T]$,  $B=1/c_{\alpha}$, $C=c_X, \delta^\pm=-D \pm \sqrt{R}$, 
with
$D = -c_h /(2c_{\alpha})$ and
$R=D^2+BC$. 

Subsequently:
\[
\mathbb{E}[X_t]=x_0e^{-\int_0^t\frac{\bar\eta(s)}{c_\alpha}ds}.
\]
From the FBSDE system for $(X_t, Y_t, Z_t)$ and centering $X_t$ and $Y_t$, one gets:
\begin{align*}
    Y_t&=\eta(t)X_t+\psi(t),\\
   \eta(t)&=\frac{c_\alpha c_g}{c_\alpha +c_g(T-t)},\\
    Z_t&=\sigma \eta(t),\\
    \psi(t)&=\big(\bar\eta(t)-\eta(t)\big)\mathbb{E}[X_t].
\end{align*}
Finally, we recall that the optimal control is given by:
\[
\hat\alpha_t=-\frac{1}{c_\alpha}Y_t=-\frac{1}{c_\alpha}\left(\eta(t)X_t+\psi(t)\right).
\]
Assuming that $X_0$ is ${\mathcal N}(x_0,\sigma_0^2)$-distributed and independent of $W$, $X_t$ is normally-distributed with mean given above by $\mathbb{E}[X_t]=x_0e^{-\int_0^t\frac{\bar\eta(s)}{c_\alpha}ds}$ and variance easily computed from
\[
dX_t=-\frac{1}{c_\alpha}\left(\eta(t)X_t+\psi(t)\right)dt+\sigma dW_t,
\]
to obtain:
\begin{align*}
     Var (X_t)&=\sigma_0^2e^{-\frac{2}{c_\alpha}\int_0^t\eta(s)ds}
  + \sigma^2\int_0^t
    e^{-\frac{2}{c_\alpha}\int_s^t\eta(s')ds'}ds.
    \end{align*}

\end{document}